\documentclass[11pt,leqno]{amsart}
\usepackage[T1]{fontenc}
\usepackage{lmodern}
\usepackage{amsmath,amssymb,amsthm,mathtools}
\usepackage[hmargin=1.1in,vmargin=1.03in]{geometry}
\usepackage{microtype}
\usepackage[hidelinks]{hyperref}
\allowdisplaybreaks[2]

\newtheorem{theorem}{Theorem}[section]
\newtheorem{corollary}[theorem]{Corollary}
\newtheorem{lemma}[theorem]{Lemma}
\newtheorem{proposition}[theorem]{Proposition}
\theoremstyle{definition}

\newtheorem{example}[theorem]{Example}
\theoremstyle{remark}
\newtheorem{remark}[theorem]{Remark}

\newcommand{\R}{\mathbb R}
\newcommand{\Ric}{\operatorname{Ric}}
\newcommand{\Hess}{\operatorname{Hess}}

\newcommand{\vol}{\operatorname{vol}}

\newcommand{\sech}{\operatorname{sech}}
\newcommand{\interior}{\mathring M}
\newcommand{\ghat}{\widehat g}
\newcommand{\gstar}{g_*}

\title[Warped splitting under spectral Ricci curvature]
{Warped Product Splitting under Spectral Ricci Curvature}
\author{Yukai Sun}
\address{School of Mathematics and Statistics, Henan University, Kaifeng
475004, P. R. China, and Center for Applied Mathematics of Henan Province,
Henan University, Zhengzhou 450046, P. R. China}
\email{sunyukai@henu.edu.cn}
\author{Changliang Wang}
\address{School of Mathematical Sciences and Institute for Advanced Study,
Key Laboratory of Intelligent Computing and Applications (Ministry of
Education), Tongji University, Shanghai 200092, China}
\email{wangchl@tongji.edu.cn}
\author{Guangrui Zhu}
\address{School of Mathematical Sciences, East China Normal University, Shanghai, 200241, P. R. China}
\email{52275500052@stu.ecnu.edu.cn}
\date{}

\begin{document}
\begin{abstract}
For $n\geq3$, we prove global warped product splitting at the spectral
Ricci coefficient $(n-1)/(n-2)$, assuming completeness of
$u^{-2(n-3)/(n-2)}g$ for a positive supersolution $u$.
On a manifold without boundary, a weighted minimizing line for this
supersolution yields a splitting with complete fiber, and the given line
is radial. If the manifold has at least two ends, such a line exists
and the fiber is compact. We also prove a splitting theorem for compact
connected boundary under a weighted mean curvature inequality.
The additional completeness is automatic in dimension three.
We also give a sufficient growth criterion and higher-dimensional
counterexamples showing that this completeness assumption cannot be omitted.
As an application, we classify orientable scalar flat three-manifolds
under the corresponding spectral conditions. The nonflat models are
positive mass spatial Schwarzschild doubles in the two-ended case and
their restrictions to either side of a coordinate sphere in the compact boundary case.
\end{abstract}

\maketitle

\section{Introduction and main results}\label{sec:intro}

The Cheeger--Gromoll splitting theorem
\cite{CheegerGromoll71} is a fundamental rigidity result in Riemannian
geometry, central to the study of the global geometry and topology of
manifolds with nonnegative Ricci curvature. In recent years, several
extensions have replaced the pointwise Ricci curvature assumption by
a spectral condition, leading to new splitting results for manifolds
with and without boundary.

To recall the spectral Ricci conditions, let $(M^n,g)$ be a Riemannian manifold and use
the sign convention $\Delta=\operatorname{tr}\Hess$. For $\gamma\geq0$,
write
\begin{equation}\label{eq:rho-def}
 \rho(x):=\min\{\Ric_g(X,X):X\in T_xM,\ |X|_g=1\},\qquad
 L_\gamma :=-\gamma\Delta+\rho.
\end{equation}
When $M$ has no boundary, we say that $(M,g)$ has
\emph{nonnegative spectral $\gamma$-Ricci curvature} if
\begin{equation}\label{eq:quadratic}
 Q_\gamma(\varphi):=\int_M
 \bigl(\gamma|\nabla\varphi|^2+\rho\varphi^2\bigr)\,dV_g\geq0
 \quad\text{for every }\varphi\in C_c^\infty(M).
\end{equation}
Equivalently, $\lambda_1(L_\gamma)\geq0$, where
\[
\lambda_1(L_\gamma) := \inf \left\{ Q_\gamma(\varphi)/\int_M\varphi^2\,dV_g \mid 0\neq\varphi\in C_c^\infty(M) \right\}.
\]
This infimum need not be an eigenvalue of $L_{\gamma}$.
For $\gamma=0$, the condition reduces to $\Ric_g\geq0$, whereas for
$\gamma>0$ it allows the least Ricci eigenvalue to be negative.
The existence of a positive $C^2$ function $u$ satisfying
$L_\gamma u\geq0$ implies \eqref{eq:quadratic} for every test function
$\varphi$. This implication follows from the ground state identity
recalled in Section~\ref{sec:distance} and is part of the
Allegretto--Piepenbrink principle;
cf.~\cite{fischer-colbrie-structure-1980}.

For manifolds without boundary, Antonelli--Pozzetta--Xu \cite{APX24}
and Catino--Mari--Mastrolia--Roncoroni \cite{CMMR25} proved that a
complete connected manifold with at least two ends and
$Q_\gamma\geq0$ for some $0\leq\gamma<4/(n-1)$ splits as a
Riemannian product $\R\times\Sigma$, with $\Sigma$ compact and of
nonnegative Ricci curvature. Hong--Wang
\cite[Theorem~1.3]{HW26a} obtained a weighted line formulation in the
same range and gave a Busemann function proof. In their formulation,
a chosen positive supersolution determines the weighted length
$\int u^\gamma\,ds_g$, and the existence of a line minimizing this
length replaces the assumption of at least two ends.

For manifolds with boundary, Hong--Wang \cite{HW26b} developed a
weighted length approach using free boundary minimizing curves.
Antonelli--Li--Sweeney \cite[Theorem~1]{ALS26} considered the case
of disconnected boundary. They proved that, for $n\geq3$ and
$0\leq\gamma<(n-1)/(n-2)$, a complete connected manifold with at
least two boundary components, at least one of them compact, splits
as $[0,L]\times\Sigma$, with $0<L<\infty$ and $\Sigma$ compact,
provided a positive $C^{2,\alpha}$ function $u$ satisfies
$L_\gamma u\geq0$ and
$H_{\mathrm{out}}+\gamma\nu_{\mathrm{out}}(\log u)\geq0$.
Here $\nu_{\mathrm{out}}$ is the outward unit normal and
$H_{\mathrm{out}}$ is the corresponding mean curvature.
Their examples show that this coefficient range is sharp for the
disconnected boundary product splitting problem.

These results motivate the study of the geometry that persists at
$\gamma = (n-1)/(n-2)$ when a Riemannian product splitting
need no longer hold. In this paper, we mainly establish two warped product
splitting theorems at this coefficient, assuming that the conformal
metric $u^{-2(n-3)/(n-2)}g$ associated with the specified positive
supersolution $u$ is complete. The first theorem concerns manifolds
without boundary that admit a weighted minimizing line; the fiber
is complete and may be noncompact. The second concerns noncompact
manifolds with compact connected boundary satisfying a weighted
mean curvature inequality; the boundary itself becomes the fiber.
In both cases, the splitting is global on $M$, and the supersolution
determines the warping function.
As a principal application, scalar flatness in dimension three leads
to a complete classification whose nonflat models are spatial
Schwarzschild doubles and their regions bounded by coordinate spheres.

Throughout the paper, assume $n\geq3$ and set
\begin{equation}\label{eq:parameters}
 m:=n-1,\qquad
 \beta:=\frac{n-1}{n-2},\qquad
 \theta:=\frac{n-3}{n-2}.
\end{equation}
For a positive function $u$, distinguish the two conformal metrics
\begin{equation}\label{eq:two-metrics}
 \ghat:=u^{2\beta}g,\qquad \gstar:=u^{-2\theta}g.
\end{equation}
We write $\hat d$ for the intrinsic length distance of $\ghat$.
The additional completeness hypothesis in our theorems concerns
$\gstar$; the weighted metric $\ghat$ may be incomplete.
Following \cite[Definition~2.1]{HW26a}, a \emph{weighted minimizing
line} for $u$ is a curve $\sigma:\R\to M$, parametrized by
$g$-arclength, such that
\[
 \hat d(\sigma(s),\sigma(t))
 =\int_s^t u(\sigma(r))^\beta\,dr
 \qquad(s<t).
\]
Its total $\ghat$-length may be finite. Our first main theorem shows
that the existence of such a line forces a warped product splitting.

\begin{theorem}[Splitting along a weighted minimizing line]\label{thm:closed}
Let $(M^n,g)$ be a smooth complete connected manifold without boundary,
$n\geq3$. Suppose a positive function
$u\in C^{2,\alpha}_{\mathrm{loc}}(M)$, $0<\alpha<1$, satisfies
\begin{equation}\label{eq:main-curvature}
 -\beta\Delta u+\rho u\geq0,
\end{equation}
and $\gstar$ is complete. If $\sigma$ is a weighted minimizing line for
$u$, then there are a complete connected smooth $(n-1)$-manifold
$(\Sigma,h)$ without boundary, a positive smooth function $a$ on $\R$,
and a constant $C>0$ such that
\begin{equation}\label{eq:closed-splitting}
 (M,g)\cong(\R\times\Sigma,dt^2+a(t)^2h),\qquad
 u(t,x)=Ca(t)^{-(n-2)}.
\end{equation}
In these coordinates $\sigma(t)=(t,x_0)$ for some $x_0\in\Sigma$. Moreover,
\begin{align}
 \rho&=\Ric_g(\partial_t,\partial_t)
       =-(n-1)\frac{a''}{a}
       =\beta\frac{\Delta u}{u},\label{eq:radial-equality}\\
 \Ric_h&\geq(n-2)\big((a')^2-aa''\big)h
       \quad\text{for every }t\in\R.\label{eq:fiber-bound}
\end{align}
In particular, $L_\beta u=0$.
\end{theorem}

When $M$ has at least two ends, Lemma~\ref{lem:curves} below provides a
weighted minimizing line for the specified supersolution. Applying
Theorem~\ref{thm:closed} then gives the following consequence, the
end assumption also forces the fiber to be compact.

\begin{corollary}[Two-ended splitting]\label{cor:two-ends}
Let $(M^n,g)$ be smooth, complete, connected and without boundary,
$n\geq3$, and suppose $M$ has at least two ends. If a positive
$u\in C^{2,\alpha}_{\mathrm{loc}}(M)$, $0<\alpha<1$, satisfies
\eqref{eq:main-curvature} and $\gstar$ is complete, then the splitting
\eqref{eq:closed-splitting} and the curvature conclusions of
Theorem~\ref{thm:closed} hold. The fiber $\Sigma$ is compact, and $M$
has exactly two ends.
\end{corollary}

The weighted line formulation also includes one-ended manifolds with
noncompact fibers, as illustrated by Example~\ref{ex:hyperbolic} below.
We now turn to a one-sided counterpart of this splitting. For a
noncompact manifold with compact connected boundary, the boundary
provides the initial hypersurface, while a weighted mean curvature
inequality supplies the comparison needed to obtain rigidity.
This leads to our second main theorem.

\begin{theorem}[Compact-boundary splitting]\label{thm:boundary}
Let $(M^n,g)$ be a smooth complete connected noncompact manifold,
$n\geq3$, whose entire boundary $B=\partial M$ is nonempty, compact,
and connected. Suppose $u>0$ is $C^{2,\alpha}$ up to the boundary locally,
satisfies \eqref{eq:main-curvature}, and $\gstar$ is complete.
Let $\nu$ be the inward unit normal and use
$H=\operatorname{tr}_B\big(g(\nabla_{\cdot}\nu,\cdot)\big)$.
If
\begin{equation}\label{eq:boundary-q}
 q_B:=H+\beta\nu(\log u)\leq0\quad\text{on }B,
\end{equation}
then $u|_B=C$ is constant and
\begin{equation}\label{eq:boundary-splitting}
 (M,g)\cong([0,\infty)\times B,dt^2+a(t)^2h),\qquad
 h=g|_B,\quad a(0)=1,\quad u(t,x)=Ca(t)^{-(n-2)} .
\end{equation}
The identities \eqref{eq:radial-equality} and the inequality
\eqref{eq:fiber-bound} hold for every $t\geq0$, and $q_B=0$.
\end{theorem}

Nonconstant warping occurs in both theorems.
Section~\ref{sec:models} gives examples satisfying all their
hypotheses in every dimension $n\geq3$.
The conclusions take a particularly simple form in dimension three,
where $\beta=4/(n-1)=2$ and $\theta=0$, so that $\gstar=g$.
The additional conformal completeness assumption is therefore
automatic, and the boundaryless result can be stated directly in
terms of $Q_2\geq0$. The following corollary also records the
rigidity of positive supersolutions in the two-ended case.

\begin{corollary}[Dimension three]\label{cor:three}
A complete connected boundaryless three-manifold with at least two ends
and $Q_2\geq0$ splits as in Corollary~\ref{cor:two-ends}. In these same coordinates
every positive $C^2$ supersolution of $L_2$ is $C/a$, with $C>0$ depending on
the supersolution. For a complete connected noncompact three-manifold
with nonempty compact connected boundary, let $u>0$ be locally
$C^{2,\alpha}$ up to the boundary and let $\nu$ be the inward unit normal.
Then $\rho\geq2\Delta u/u$ and $H+2\nu(\log u)\leq0$ imply
\eqref{eq:boundary-splitting}.
No additional conformal completeness condition is needed.
\end{corollary}

For the boundaryless two-ended problem, this corollary describes
rigidity at the critical coefficient $\gamma=2$: below it one obtains
a Riemannian product \cite{APX24,CMMR25}, whereas at the endpoint
nonconstant warping is possible.

Under the additional assumption of scalar flatness, the warping
function satisfies a first integral whose constant becomes the
Schwarzschild mass after normalizing the spherical fiber.
This yields the following classification, which
we regard as a principal geometric application of the splitting
theorems. For $\mathfrak m>0$, let $\mathcal S_{\mathfrak m}$ denote
the spatial Schwarzschild double, obtained by smoothly doubling
\[
 g_{\mathfrak m}
 =\frac{dr^2}{1-2\mathfrak m/r}+r^2g_{S^2(1)},
 \qquad r>2\mathfrak m,
\]
across its minimal sphere $r=2\mathfrak m$.

\begin{corollary}[Scalar flat three-dimensional classification]
\label{cor:scalar-flat-classification}
Let $(M^3,g)$ be a smooth, complete, connected, orientable
Riemannian manifold with $R_g=0$.
\begin{enumerate}
\item[(i)]
If $M$ has no boundary and has at least two ends, then $Q_2\geq0$
if and only if $(M,g)$ is isometric either to
$(\R\times T^2,dt^2+h)$, where $h$ is flat, or to
$\mathcal S_{\mathfrak m}$ for some $\mathfrak m>0$.

\item[(ii)]
Suppose $M$ is noncompact and its entire boundary $B=\partial M$
is nonempty, compact, and connected. Let $H$ be the sum of the
principal curvatures with respect to the inward unit normal.
Then the Robin quadratic form satisfies
\[
 Q_B(\varphi):=
 \int_M\bigl(2|\nabla\varphi|^2+\rho\varphi^2\bigr)\,dV_g
 -\int_B H\varphi^2\,dA_g\geq0
\]
for every smooth compactly supported function up to $B$, with
arbitrary boundary trace, if and only if $(M,g)$ is isometric to
one of the following: a flat product $[0,\infty)\times T^2$;
the closed exterior of a round ball in Euclidean three-space;
or the closure of either component obtained by cutting
$\mathcal S_{\mathfrak m}$, $\mathfrak m>0$, along a coordinate
sphere, including its minimal sphere.
\end{enumerate}
In case~\textup{(ii)}, if $Q_B\geq0$ and $B$ is a minimal
two-sphere, then $(M,g)$ is the Schwarzschild exterior bounded
by its horizon, with
\[
 \mathfrak m=\sqrt{\frac{|B|_g}{16\pi}}.
\]
\end{corollary}

This result gives an intrinsic spectral characterization of the
Schwarzschild geometries. In the nonflat cases, spherical symmetry
and asymptotic flatness follow from the hypotheses, without an
asymptotic model or a static potential being prescribed.
In particular, the minimal spherical boundary case determines
both the complete geometry and its mass. These conclusions
complement classifications based on the static vacuum equations;
see, for example, \cite{ReirisI,ReirisII}.
The proofs are given in Propositions~\ref{prop:schwarz-closed}
and~\ref{prop:schwarz-boundary} in Section~\ref{sec:applications}.

In dimensions $n>3$, the two spectral coefficients appearing in
the preceding discussion are distinct:
\begin{equation}\label{eq:threshold-gap}
 \frac{n-1}{n-2}-\frac4{n-1}
 =\frac{(n-3)^2}{(n-1)(n-2)}>0
 \qquad(n>3).
\end{equation}
Thus our higher-dimensional results at $\beta$ extend beyond the
boundaryless product splitting range by imposing the additional
completeness of $\gstar$. The following remark records why this
hypothesis is essential to the statements.

\begin{remark}[The conformal completeness hypothesis]
\label{rem:main-completeness}
For every $n\geq4$, the completeness of
$g_*=u^{-2(n-3)/(n-2)}g$ cannot be omitted from
Theorems~\ref{thm:closed} and~\ref{thm:boundary}.
Proposition~\ref{prop:beta-counterexample} below gives complete manifolds satisfying
all the other hypotheses, including the boundary inequality in the
compact boundary case, which admit no local warped product representation
with one-dimensional base. In these examples $g_*$ is incomplete.
This condition is not necessary for every individual warped product;
see Remark~\ref{rem:completeness-nonremovable}. In particular, the
splitting conclusion itself does not force $\gstar$ to be complete.
\end{remark}

Failure of warped product rigidity under the spectral condition
alone already occurs at the smaller coefficient $c=4/(n-1)$.
For every $n\geq4$, Proposition~\ref{prop:critical-all} below constructs
a complete two-ended metric on $\R\times S^{n-1}$ with $Q_c\geq0$
that admits no local warped product splitting with one-dimensional
base near a specified point. Since
\[
 Q_\gamma(\varphi)
 =Q_c(\varphi)+(\gamma-c)\int_M|\nabla\varphi|^2\,dV_g,
\]
the same example satisfies $Q_\gamma\geq0$ for every $\gamma\geq c$.
In particular, spectral nonnegativity alone cannot replace the
conformal completeness hypothesis at $\gamma=\beta$.
On the other hand, this completeness follows from simple upper
growth bounds on $u$, as the next criterion shows.

\begin{proposition}[Integral growth criterion]\label{prop:growth}
Let $n>3$, let $(M,g)$ be complete, with or without boundary, fix
$o\in M$, and set $r(x)=d_g(o,x)$. If a positive nondecreasing function
$A:[0,\infty)\to(0,\infty)$ satisfies
\begin{equation}\label{eq:A-condition}
 u(x)\leq A(r(x)),\qquad
 \int_1^\infty A(s)^{-\frac{n-3}{n-2}}\,ds=\infty ,
\end{equation}
then $\gstar$ is complete.
In particular, an upper bound $u\leq C$ suffices, as does
\begin{equation}\label{eq:polynomial-condition}
 u(x)\leq C(1+r(x))^{\frac{n-2}{n-3}} .
\end{equation}
\end{proposition}

\begin{proof}
Define $F(R)=\int_0^R A(s)^{-\theta}\,ds$.
The radial function $r$ is $1$-Lipschitz for the intrinsic metric of $g$,
also when a boundary is present. For every piecewise smooth curve $c$
from $o$ to $x$, absolute continuity and the chain rule give
\begin{align*}
 L_{\gstar}(c)
 &=\int u(c)^{-\theta}|\dot c|_g
 \geq\int A(r(c))^{-\theta}|\dot c|_g\\
 &\geq\int A(r(c))^{-\theta}|(r\circ c)'|
 \geq F(r(x)).
\end{align*}
Note that $A$ need not be continuous. The integrand $A^{-\theta}$
is locally bounded and has only countably many discontinuities.
Away from these values $F'=A^{-\theta}$; on the preimage of each
exceptional value the absolutely continuous function $r\circ c$ has
derivative zero almost everywhere. Since $F$ is locally Lipschitz,
this gives the stated chain rule also on those preimages. Taking the infimum yields
$d_{\gstar}(o,x)\geq F(r(x))$. Since $F(R)\to\infty$, every bounded
$\gstar$-ball lies in a bounded $g$-ball. Closed $g$-balls are compact
because $g$ is complete, and the two metrics induce the same topology;
hence every closed $\gstar$-ball, being a closed subset of a compact
$g$-ball, is compact. Thus $\gstar$ is proper and therefore complete. For \eqref{eq:polynomial-condition}
one can take $A(s)=C(1+s)^{1/\theta}$, for which the integral is a
constant multiple of $\int_1^\infty(1+s)^{-1}\,ds$.
\end{proof}

We briefly describe the common mechanism behind Theorems \ref{thm:closed} and \ref{thm:boundary}. The proof begins with weighted second variation and
distance comparison, using complete approximations of the weighted
metric when necessary. Completeness of $\gstar$ controls the error
terms in the comparison and allows us to construct a smooth
weighted Busemann function $b$ satisfying
\[
 \Delta_g b=0,\qquad |\nabla b|_g=u^\beta.
\]
In the boundaryless case, $b$ is obtained from the two directions of
the given weighted line. In the boundary case, the weighted
mean curvature inequality gives a comparison for the distance from
the boundary $\partial M$, and a minimizing ray forces that distance to be harmonic.
These constructions are carried out in
Proposition~\ref{prop:harmonic}; see also Remark~\ref{rem:busemann}.

Once this harmonic function is available, the Bochner formula and
the equality case of the refined Kato inequality determine the
geometry of its level sets and yield the warped product structure,
as in \cite[Lemma~4.1]{Li-Wang-ASENS-weighted}.
A global flow argument then identifies the whole manifold with the
resulting warped product. The conformal curvature identity also
relates our setting to Wylie's weighted warped splitting theorem
\cite{Wylie17}. Since $\ghat$ may be incomplete, however, the
distance comparison required here is established directly.

As an application, we classify complete orientable scalar flat
three-manifolds in the two-ended and compact connected boundary
settings under the corresponding spectral conditions.
The nonflat two-ended models are positive mass spatial Schwarzschild
doubles, and the nonflat compact boundary models are the regions
obtained by cutting such doubles along coordinate spheres.
For the boundary classification, the supersolution and weighted
mean curvature conditions can be expressed through a Robin quadratic
form; see Proposition~\ref{prop:schwarz-boundary} below.

The paper is organized as follows.
Section~\ref{sec:distance} develops the positive supersolution and
conformal identities and establishes the weighted distance
comparisons. Section~\ref{sec:busemann} constructs the weighted
Busemann functions and proves their harmonicity and gradient
identities. Section~\ref{sec:rigidity} analyzes the equality case
and the global normal flow, completing the proofs of the main
splitting theorems. Section~\ref{sec:models} presents nonconstant
warped examples, and Section~\ref{sec:counterexamples} gives the
counterexamples and examines the role of conformal completeness
and the criticality of the associated operators.
Finally, Section~\ref{sec:applications} proves the scalar flat
three-dimensional classifications.

\medskip
\noindent\textbf{Acknowledgements.}
We thank Dr.~Yifan Xiao for discussing a counterexample to warped product
splitting at the coefficient $(n-1)/(n-2)$.
Y. S. is partially supported by the National Natural Science
Foundation of China (No. 12601093). C. W. is partially supported by the National Natural Science
Foundation of China (No. 12671068).

\subsection*{Declaration of generative AI use}

An AI assistant was used for language editing and LaTeX preparation of this
manuscript. All mathematical results and proofs are due to the authors, who take
full responsibility for the content of the article.

\section{Positive supersolutions and weighted distance}\label{sec:distance}

This section prepares the analytic tools used throughout the paper.
We first record the ground-state identity and the existence of positive
supersolutions on exhausted domains. We then establish the conformal
tensor identity relating the Bakry--\'Emery Ricci tensor of the weighted
metric $\ghat$ to the spectral Ricci bound. Combining this identity with
the index form of $\ghat$-geodesics yields a weighted distance comparison,
which is the key input for the Busemann construction in the next section.

For later use, the ground-state identity is
\begin{equation}\label{eq:ground-state}
 Q_\gamma(\varphi)
 =\gamma\int_Mu^2\left|\nabla\left(\frac{\varphi}{u}\right)\right|^2dV_g
  +\int_M\frac{L_\gamma u}{u}\varphi^2\,dV_g .
\end{equation}
We derive \eqref{eq:ground-state}. Put $\psi=\varphi/u$ and
\[
 \nabla\varphi=\psi\nabla u+u\nabla\psi,\qquad
 |\nabla\varphi|^2=\psi^2|\nabla u|^2+2u\psi\langle\nabla u,\nabla\psi\rangle
                  +u^2|\nabla\psi|^2 .
\]
Since $2u\psi\langle\nabla u,\nabla\psi\rangle
       =\langle\nabla u,\nabla(u\psi^2)\rangle-\psi^2|\nabla u|^2$, this gives
\[
 \gamma|\nabla\varphi|^2+\rho\varphi^2
 =\gamma u^2|\nabla\psi|^2+\gamma\langle\nabla u,\nabla(u\psi^2)\rangle
  +\rho u^2\psi^2 .
\]
Integrate over $M$. The middle term is, by the divergence theorem,
\[
 \gamma\int_M\langle\nabla u,\nabla(u\psi^2)\rangle\,dV_g
 =-\gamma\int_M(\Delta u)\,u\psi^2\,dV_g
 =-\gamma\int_M\frac{\Delta u}{u}\,\varphi^2\,dV_g ,
\]
with no boundary term because $\varphi$ has compact support. Hence
\[
 Q_\gamma(\varphi)
 =\gamma\int_Mu^2|\nabla\psi|^2\,dV_g
  +\int_M\Big(\rho-\gamma\frac{\Delta u}{u}\Big)\varphi^2\,dV_g
 =\gamma\int_Mu^2\left|\nabla\left(\frac{\varphi}{u}\right)\right|^2dV_g
  +\int_M\frac{L_\gamma u}{u}\varphi^2\,dV_g ,
\]
which is \eqref{eq:ground-state}. On a manifold with boundary, the
same identity holds for $\varphi\in C_c^\infty(\interior)$.
For test functions with arbitrary boundary trace we retain the
boundary term in Proposition~\ref{prop:schwarz-boundary}.

We recall why $Q_\gamma\geq0$, $\gamma>0$, gives a positive zero
solution on a connected, noncompact manifold without boundary. This is
the case needed in Corollary~\ref{cor:three}. Fix a smooth connected
exhaustion $\Omega_j\Subset\Omega_{j+1}$ and $o\in\Omega_1$.

\smallskip\noindent
\emph{Step 1: strict Dirichlet positivity.}
Let $\lambda_1(\Omega_j)$ be the first Dirichlet eigenvalue of $L_\gamma$
on $\Omega_j$. Then $\lambda_1(\Omega_j)>0$. Indeed, if
$\lambda_1(\Omega_j)=0$, its positive first eigenfunction, extended by
zero to $\Omega_{j+1}$, would minimize the nonnegative form
$Q_\gamma$ with value zero, hence be a weak solution on $\Omega_{j+1}$.
The strong maximum principle would then force it to be positive in the
interior, contradicting its vanishing on
$\Omega_{j+1}\setminus\Omega_j\neq\varnothing$.

\smallskip\noindent
\emph{Step 2: a positive solution on each $\Omega_j$.}
Since $\lambda_1(\Omega_j)>0$, the Dirichlet operator is invertible, so
there is a unique $u_j$ with
\[
 L_\gamma u_j=0\ \text{ in }\Omega_j,\qquad u_j=1\ \text{ on }\partial\Omega_j .
\]
The negative part $u_j^-$ has zero trace. Testing the equation against
$u_j^-$ gives $Q_\gamma(u_j^-)=0$, so $u_j^-=0$ by Step 1. Thus
$u_j\geq0$, and Harnack's inequality with the nonzero boundary data gives
$u_j>0$.

\smallskip\noindent
\emph{Step 3: passage to the limit.}
Normalize $u_j(o)=1$. On each precompact connected subdomain
containing $o$, interior Harnack bounds and local elliptic estimates give
a subsequential $C^2_{\rm loc}$ limit. A diagonal argument yields
$u\in C^2_{\rm loc}(M)$ with
\[
 L_\gamma u=0\ \text{ on }M,\qquad u(o)=1,\qquad u>0\ \text{(Harnack)} .
\]
Since $\rho$ is locally Lipschitz, elliptic regularity upgrades $u$ to
$C^{2,\alpha}_{\mathrm{loc}}$ for every $0<\alpha<1$.

\smallskip
This proves the analytic assertion needed for Corollary~\ref{cor:three}.
The distance computations below use only two derivatives of $u$, so they
apply at this regularity. After splitting, the identity
$u=|\nabla b|^{1/\beta}$, where $b$ is the harmonic function constructed in
Proposition~\ref{prop:harmonic}, also proves its smoothness.

The distance between two points and, when a boundary $B$ is present,
the distance to $B$ are given by
\begin{equation}\label{eq:D-def}
 \hat{d}(p,q)=\inf_c\int_c u^\beta\,ds_g,\qquad
 \hat{d}(B,q)=\inf_{p\in B}\hat{d}(p,q).
\end{equation}
Even if $\ghat$ is incomplete, $\hat{d}$ is a distance inducing the manifold
topology and is locally Lipschitz in each variable.
For example, a sufficiently small fixed coordinate neighborhood has
$u$ bounded above and below by positive constants; exiting a smaller
neighborhood costs a positive amount of weighted length.
These observations establish positivity of $\hat{d}$ and the local assertions.

\subsection{The conformal tensor identity}

Set $w=\log u$, $\phi=\beta w$, and $f=mw$.
A direct conformal change computation gives
\begin{equation}\label{eq:conformal-tensor}
 \Ric_{\ghat}+\Hess_{\ghat}f+\frac1m\,df\otimes df
 =
 \Ric_g-\beta\frac{\Delta_gu}{u}g
 +\beta\big(|dw|_g^2g-dw\otimes dw\big)\geq0 .
\end{equation}
We record the formulas used:
\begin{align*}
 \Ric_{e^{2\phi}g}
 &=\Ric_g-(n-2)(\Hess_g\phi-d\phi\otimes d\phi)
       -(\Delta_g\phi+(n-2)|d\phi|^2)g,\\
 \Hess_{e^{2\phi}g}f
 &=\Hess_gf-df\otimes d\phi-d\phi\otimes df
       +\langle df,d\phi\rangle_g g .
\end{align*}
Adding $df\otimes df/m$, using $(n-2)\beta=m$, and then
$\Delta w=\Delta u/u-|dw|^2$ proves \eqref{eq:conformal-tensor}.
The associated drift Laplacian satisfies
\begin{equation}\label{eq:laplace-transform}
 \Delta_g F=e^{2\phi}\left(
   \Delta_{\ghat}F-\langle\nabla_{\ghat}f,\nabla_{\ghat}F\rangle_{\ghat}
                         \right).
\end{equation}
The tensor on the left of \eqref{eq:conformal-tensor} is the
$N=1$ Bakry--\'Emery tensor of $(\ghat,f)$.
The length parameter appearing in weighted comparison is
\begin{equation}\label{eq:affine-length}
 e^{-2f/m}\,d\tau=u^{-2}u^\beta\,ds_g
                  =u^{-\theta}\,ds_g=ds_{\gstar},
\end{equation}
where $\tau$ is $\ghat$-arclength. This explains the relation to
the $f$-completeness condition in \cite[Section~6]{Wylie17}.
We assume neither that $\ghat$ is complete nor that $f$ is bounded
above, so we use the following finite segment argument instead of
invoking that splitting theorem. To illustrate the distinction between the two lengths in
\eqref{eq:affine-length}, let $(\mathbb T^2,h)$ be a flat torus and set
$$
  M=\mathbb R\times\mathbb T^2,\qquad
  g=dt^2+e^{2t}h,\qquad u=e^{-t}.
$$
The metric $g$ is complete and has two ends. The warped product formulas
give $\Ric_g=-2g$ and $\Delta_g u=u''+2u'=-u$, so $L_2u=0$.
In dimension three $\beta=2$ and $\theta=0$, hence $\gstar=g$ is
complete, while $f=2\log u=-2t$ is unbounded above as $t\to-\infty$.
At the other end, fix $x_0\in\mathbb T^2$ and consider the radial ray
$\sigma(t)=(t,x_0)$ for $t\geq0$. Along this ray,
$$
  L_{\ghat}(\sigma)=\int_0^\infty e^{-2t}\,dt=\frac12.
$$
Thus $\ghat$ is incomplete even though the weighted length entering
the comparison is infinite along this ray.

\subsection{Index comparison on finite segments}

An inequality $\Delta F\leq c$ in the upper support sense means that
at each point and for each $\varepsilon>0$ there is a $C^2$ upper
support agreeing there and having Laplacian at most $c+\varepsilon$.
Indeed, if a $C^2$ function $\psi$ touches $F$ from below at $q$ and
$S_\varepsilon$ is such an upper support, then
$S_\varepsilon-\psi$ has a local minimum at $q$, whence
$\Delta_g\psi(q)\leq\Delta_g S_\varepsilon(q)\leq c(q)+\varepsilon$.
Letting $\varepsilon$ tend to zero gives the viscosity inequality.
For the smooth linear operator $\Delta_g$ and continuous right hand
side this implies the distributional inequality
\cite{Ishii1995weak}, the reversed statement follows by changing sign.
These distributional bounds pass to locally uniform limits by testing
against nonnegative compactly supported smooth functions.

\begin{lemma}[Point-distance comparison]\label{lem:point-comparison}
Let $(M^n,g)$ be a smooth complete connected manifold, possibly with boundary, $n\geq3$. Suppose a positive function
$u\in C^{2,\alpha}_{\mathrm{loc}}(M)$, $0<\alpha<1$, satisfies
 $-\beta\Delta u+\rho u\geq0$. If an interior
$\hat{d}$-minimizing segment $\sigma$ joins $p$ to $q$, then
\begin{equation}\label{eq:segment-comparison}
 \Delta_{g,q}\hat{d}(p,q)
 \leq
 \frac{m\,u(q)^{2/(m-1)}}{\displaystyle\int_\sigma u^{-\theta}\,ds_g}
\end{equation}
in the upper support sense. If $g$ and $\gstar$ are complete and
$M$ has no boundary, then for every $p\ne q$
\begin{equation}\label{eq:global-comparison}
 \Delta_{g,q}\hat{d}(p,q)
 \leq\frac{m\,u(q)^{2/(m-1)}}{d_{\gstar}(p,q)}
\end{equation}
in the same sense, whether or not the infimum defining $\hat{d}(p,q)$
is attained.
\end{lemma}

\begin{proof}
First suppose a minimizing segment exists, and parametrize it by
$\ghat$-arclength $\tau\in[0,L]$. Write $T=\dot\sigma$.
Let $E_1,\ldots,E_m$ be a parallel orthonormal normal frame.
The second variation with trial fields $\xi E_i$, where
$\xi(0)=0$ and $\xi(L)=1$, gives an upper support $F$ at $q$ with
$\nabla_{\ghat}F(q)=T(L)$ and
\begin{equation}\label{eq:index-start}
 e^{-2\phi(q)}\Delta_gF(q)
 \leq\int_0^L\big(m(\xi')^2-\xi^2\Ric_{\ghat}(T,T)\big)\,d\tau
         -f'(L).
\end{equation}
Here is a simultaneous variation giving this support, also at a cut
point. Let $k$ be a smooth auxiliary metric, put $V_i=\xi E_i$
for $i\leq m$ and $V_{m+1}=(\tau/L)T$, and set
\begin{equation}\label{eq:simultaneous-variation}
 c_z(\tau)=\exp^k_{\sigma(\tau)}\!\left(\sum_{i=1}^{m+1}z_iV_i(\tau)\right),
 \qquad |z|<\varepsilon_0 .
\end{equation}
Since $c_z(0)=p$ and $\partial_{z_i}c_z(L)\big|_{z=0}=V_i(L)$, the
differential of the endpoint map $\mathcal E(z)=c_z(L)$ at $z=0$
sends the standard basis to the $\ghat$-orthonormal basis
$E_1(L),\ldots,E_m(L),T(L)$ of $T_qM$. Hence $\mathcal E$ is a
diffeomorphism of a ball $B_{\varepsilon_0}(0)$ onto a neighborhood
$U$ of $q$, and we define $F$ explicitly by
\begin{equation}\label{eq:support-definition}
 F(x)=L_{\ghat}\big(c_{\mathcal E^{-1}(x)}\big)
 =\int_0^L\big|\partial_\tau c_{\mathcal E^{-1}(x)}(\tau)\big|_{\ghat}\,d\tau,
 \qquad x\in U .
\end{equation}
Every $c_z$ joins $p$ to $\mathcal E(z)$, so $F\geq \hat{d}(p,\cdot)$ on
$U$; since $\sigma$ is minimizing, $F(q)=L=\hat{d}(p,q)$. For small $z$
the speeds stay nonzero, uniformly on $[0,L]$; since $\ghat$ is
$C^{2,\alpha}$, differentiating the length integral twice in $z$
gives continuous derivatives. Thus $F$ is a $C^2$ upper support at
$q$, whether or not $q$ is a cut point of $p$.

Write $\delta_i$ for differentiation in the variation parameter $z_i$.
For the endpoint function $F$ we set
$\delta_iF=\partial_{z_i}(F\circ\mathcal E)(0)$ and
$\delta_j\delta_iF=\partial_{z_j}\partial_{z_i}(F\circ\mathcal E)(0)$.
Since $\sigma$ is a $\ghat$-geodesic, the first variation of length gives
\begin{equation}\label{eq:first-variation}
 \delta_iF
 =\big\langle T,V_i\big\rangle_{\ghat}\Big|_0^L
 =\ghat\big(T(L),V_i(L)\big),
\end{equation}
because $V_i(0)=0$ for every $i$. In particular
$dF_q=\ghat(T(L),\cdot)$.

For the second derivatives, all covariant derivatives are taken with
respect to $\ghat$; the auxiliary metric $k$ serves only to define
the variation \eqref{eq:simultaneous-variation} smoothly. Write
$V_i^\perp=V_i-\ghat(V_i,T)T$, and use the curvature convention
$\Ric_{\ghat}(T,T)=\sum_i\langle R^{\ghat}(E_i,T)T,E_i\rangle_{\ghat}$.
The first variation integrand is
$\langle\widehat\nabla_{\partial_{z_i}c_z}\partial_\tau c_z,
 \partial_\tau c_z\rangle_{\ghat}/|\partial_\tau c_z|_{\ghat}$.
Differentiating it in $z_j$ at $z=0$, where
$\partial_{z_i}c_0=V_i$, $|\partial_\tau c_0|_{\ghat}=1$, and
$\widehat\nabla_{V_i}T=\widehat\nabla_TV_i$ by torsion freeness, gives
\[
 \delta_j\delta_iF
 =\int_0^L\Big(
   \big\langle\widehat\nabla_{V_j}\widehat\nabla_{V_i}T,\,T
   \big\rangle_{\ghat}
   +\big\langle\widehat\nabla_TV_i,\widehat\nabla_TV_j\big\rangle_{\ghat}
   -\big\langle\widehat\nabla_TV_i,T\big\rangle_{\ghat}
    \big\langle\widehat\nabla_TV_j,T\big\rangle_{\ghat}
   \Big)\,d\tau .
\]
Torsion freeness gives
$\widehat\nabla_{V_j}\widehat\nabla_{V_i}T
 =\widehat\nabla_T\widehat\nabla_{V_i}V_j
  +R^{\ghat}(V_j,T)V_i$
and $\widehat\nabla_{V_i}V_j=\widehat\nabla_{V_j}V_i$, since
$\widehat\nabla_TT=0$, integration by parts yields
\[
 \int_0^L\big\langle\widehat\nabla_{V_j}\widehat\nabla_{V_i}T,\,T
 \big\rangle_{\ghat}\,d\tau
 =\big\langle\widehat\nabla_{V_j}V_i,\,T\big\rangle_{\ghat}
   \Big|_0^L
  -\int_0^L\big\langle R^{\ghat}(V_i,T)T,\,V_j\big\rangle_{\ghat}\,d\tau .
\]
Hence
\begin{align}
 \delta_j\delta_iF
 &=\int_0^L\Big(
    \big\langle\widehat\nabla_TV_i,\widehat\nabla_TV_j\big\rangle_{\ghat}
   -\big\langle\widehat\nabla_TV_i,T\big\rangle_{\ghat}
    \big\langle\widehat\nabla_TV_j,T\big\rangle_{\ghat}
   -\big\langle R^{\ghat}(V_i,T)T,V_j\big\rangle_{\ghat}
   \Big)\,d\tau \notag\\
 &\quad+\big\langle\widehat\nabla_{V_j}V_i,\,T\big\rangle_{\ghat}(L),
 \label{eq:second-variation-length}
\end{align}
the boundary term at $\tau=0$ vanishing because $c_z(0)=p$ for all
$z$. Since $\widehat\nabla_TT=0$, the splitting into normal and
tangential parts is parallel along $\sigma$, so
\[
 \big\langle\widehat\nabla_TV_i,\widehat\nabla_TV_j\big\rangle_{\ghat}
 -\big\langle\widehat\nabla_TV_i,T\big\rangle_{\ghat}
  \big\langle\widehat\nabla_TV_j,T\big\rangle_{\ghat}
 =\big\langle\widehat\nabla_TV_i^\perp,
           \widehat\nabla_TV_j^\perp\big\rangle_{\ghat},
\]
and $\langle R^{\ghat}(V_i,T)T,V_j\rangle_{\ghat}
 =\langle R^{\ghat}(V_i^\perp,T)T,V_j^\perp\rangle_{\ghat}$.
The integral in \eqref{eq:second-variation-length} is therefore the
index form $I(V_i^\perp,V_j^\perp)$; in particular the tangential
field $V_{m+1}=(\tau/L)T$ contributes zero to it.

It remains to convert $\delta_j\delta_iF$ into the Hessian of $F$.
Differentiating $F(\mathcal E(z))=L_{\ghat}(c_z)$ twice by the chain
rule,
\begin{equation}\label{eq:hess-chain-rule}
 \delta_j\delta_iF
 =\Hess_{\ghat}F\big(V_j(L),V_i(L)\big)
  +\big\langle\nabla_{\ghat}F,\,
           \widehat\nabla_{V_j}V_i\big\rangle_{\ghat}(L).
\end{equation}
By \eqref{eq:first-variation},
$\langle\nabla_{\ghat}F,\widehat\nabla_{V_j}V_i\rangle_{\ghat}(L)
 =\langle T,\widehat\nabla_{V_j}V_i\rangle_{\ghat}(L)$.
The terminal acceleration term therefore appears with the same sign in
\eqref{eq:second-variation-length} and \eqref{eq:hess-chain-rule}.
Subtracting this common term gives
\begin{equation}\label{eq:hess-equals-index}
 \Hess_{\ghat}F\big(V_j(L),V_i(L)\big)
 =I(V_i^\perp,V_j^\perp)
 =\int_0^L\Big(
   \big\langle\widehat\nabla_TV_i^\perp,\widehat\nabla_TV_j^\perp
   \big\rangle_{\ghat}
   -\big\langle R^{\ghat}(V_i^\perp,T)T,V_j^\perp\big\rangle_{\ghat}
   \Big)\,d\tau .
\end{equation}
Tracing \eqref{eq:hess-equals-index} in the terminal orthonormal
frame, the tangential direction contributes zero, and
$\widehat\nabla_TV_i=\xi'E_i$ gives
\[
 \Delta_{\ghat}F(q)
 =\sum_{i=1}^m\Hess_{\ghat}F(E_i(L),E_i(L))
 =\int_0^L\big(m(\xi')^2-\xi^2\Ric_{\ghat}(T,T)\big)\,d\tau .
\]
Since $\nabla_{\ghat}F(q)=T(L)$, the drift term in
\eqref{eq:laplace-transform} is
$\langle\nabla_{\ghat}f,\nabla_{\ghat}F\rangle=f'(L)$, and
\eqref{eq:laplace-transform} gives \eqref{eq:index-start}. Thus no
differentiability of the distance or of a choice of minimizing
segments in their endpoints is used.

By \eqref{eq:conformal-tensor}, evaluated on the unit tangent $T$,
$\Ric_{\ghat}(T,T)\geq-f''-(f')^2/m$. Substituting into
\eqref{eq:index-start} and integrating by parts,
\[
 \int_0^L\xi^2 f''\,d\tau
 =\big[\xi^2 f'\big]_0^L-2\int_0^L\xi\xi' f'\,d\tau
 =f'(L)-2\int_0^L\xi\xi' f'\,d\tau ,
\]
since $\xi(0)=0$, $\xi(L)=1$. The boundary term $f'(L)$ cancels the
$-f'(L)$ in \eqref{eq:index-start}, leaving
\begin{equation}\label{eq:index-square}
 e^{-2\phi(q)}\Delta_gF(q)
 \leq\int_0^L\Big(m(\xi')^2-2\xi\xi' f'+\frac{(f')^2}{m}\xi^2\Big)\,d\tau
 =m\int_0^L\left(\xi'-\frac{f'}m\xi\right)^2\,d\tau .
\end{equation}

Now let $s\in[0,\ell]$ be $g$-arclength and write $\tau=\tau(s)$.
Along $\sigma$ we have $d\tau=u^\beta ds$. Set
$\xi(\tau(s))=u(\sigma(s))\zeta(s)/u(q)$, where
$\zeta(0)=0$ and $\zeta(\ell)=1$.
Since $f=m\log u$, we have $\partial_\tau f=m u^{-\beta}(\partial_s u)/u$.
The product rule gives
$$
 \partial_\tau\xi-\frac{\partial_\tau f}{m}\xi
 =u^{-\beta}\frac{u}{u(q)}\,\partial_s\zeta .
$$
Using $2-\beta=\theta$ and writing $\zeta'=d\zeta/ds$, we obtain
$$
 m\int_0^L\left(\partial_\tau\xi-\frac{\partial_\tau f}{m}\xi\right)^2d\tau
 =\frac{m}{u(q)^2}\int_0^\ell u^\theta(\zeta')^2\,ds .
$$
Multiplying \eqref{eq:index-square} by $e^{2\phi(q)}=u(q)^{2\beta}$ and
using $2\beta-2=2/(m-1)$ gives
\begin{equation}\label{eq:energy-comparison}
 \Delta_gF(q)\leq
 m u(q)^{2/(m-1)}\int_0^\ell u^\theta(\zeta')^2\,ds .
\end{equation}

The weighted Dirichlet energy $\int_0^\ell u^\theta(\zeta')^2ds$ with
$\zeta(0)=0,\zeta(\ell)=1$ is minimized by the solution of
$(u^\theta\zeta')'=0$, i.e. $u^\theta\zeta'=\text{const}$. Hence
\begin{equation}\label{eq:optimal-cutoff}
 \zeta(s)=\frac{\int_0^s u^{-\theta}\,dr}{\int_0^\ell u^{-\theta}\,dr},
 \qquad
 \min\int_0^\ell u^\theta(\zeta')^2\,ds
 =\left(\int_0^\ell u^{-\theta}\,ds\right)^{-1}.
\end{equation}
Substituting into \eqref{eq:energy-comparison} proves
\eqref{eq:segment-comparison}; since the denominator is the
$\gstar$-length of $\sigma$, \eqref{eq:global-comparison} follows when a
minimizer exists.

Here are the details for a nonattained distance.
Choose nested compact sets $K_j\subset\operatorname{int}_M K_{j+1}$
exhausting $M$ and smooth cutoffs $0\leq\chi_j\leq1$ equal to one near
$K_j$ and supported in $K_{j+1}$. Put $u_j=u+1-\chi_j$ and let $\hat{d}_j$
be the intrinsic length distance for $\widehat g_j=u_j^{2\beta}g$.
For each $j$ the function $u_j$ has a positive global lower bound:
it is at least one outside $K_{j+1}$ and positive on that compact set.
Thus $\hat{d}_j$ is complete. It is a locally compact length metric, hence
proper and geodesic by the length-space Hopf--Rinow theorem.
We have $\hat{d}_j\geq \hat{d}$. For fixed endpoints any competitor with
$\ghat$-length at most $\hat{d}+\varepsilon$ has compact image, hence the
same length for $\widehat g_j$ once $j$ is large. This proves
pointwise convergence. On precompact coordinate balls (half balls
at the boundary) the triangle inequality and local upper bounds on
$u_j$ make $\hat{d}_j$ equicontinuous in both variables: nearby points
can be joined inside the chart with length bounded by a fixed multiple
of their coordinate separation. A finite cover on compact sets yields
$\hat{d}_j\to \hat{d}$ locally uniformly on $M\times M$. This argument applies
unchanged to intrinsic distances when a boundary is included.

Take $\hat{d}_j$-minimizing curves from $q$ to $p$, parametrized by
$g$-arclength. If their $g$-lengths stay bounded, completeness of $g$
and the geodesic equation on compact sets give a $\hat{d}$-minimizing limit.
Otherwise take a subsequence with lengths tending to infinity.
For any fixed $T$ the arcs of $g$-length $T$ from $q$ lie in the
compact ball $\overline B_g(q,T)$, where eventually $u_j=u$.
With $Y=\dot\sigma_j$ their reparametrized geodesic equation there is
$\nabla_Y Y=\nabla\phi-d\phi(Y)Y$, $|Y|_g=1$, by conformal change.
Compactness of the initial directions and continuous dependence for
this equation give convergence on each finite parameter interval to
a $g$-unit-speed curve $\sigma:[0,\infty)\to M$ with $\sigma(0)=q$.
Each limiting segment minimizes $\hat{d}$: its length is the limit of the
segment lengths, and $\hat{d}_j$ converges locally uniformly also at moving
endpoints. In particular,
\begin{equation}\label{eq:exact-decomposition}
 \hat{d}(q,p)=\hat{d}(q,\sigma(s))+\hat{d}(\sigma(s),p)\quad(s>0).
\end{equation}
To obtain this identity, split each approximating minimizing curve
at $g$-arclength $s$, use local convergence of that initial segment
and local uniform convergence of $\hat{d}_j$, and pass to the limit.
The ray is proper. Otherwise choose $t_k\to\infty$ with
$t_{k+1}>t_k+1$ and $\sigma(t_k)$ in some compact $C$.
The weighted lengths of $[0,t_k]$ equal $\hat{d}(q,\sigma(t_k))$, bounded
above by $\max_C \hat{d}(q,\cdot)$. Yet each interval $[t_i,t_i+1]$ lies
in the compact closed unit $g$-neighborhood $C_1$ of $C$ and contributes
at least $\min_{C_1}u^\beta>0$. Summing for $i<k$ is a contradiction.

For fixed $s$, add the constant $\hat{d}(\sigma(s),p)$ to the upper support
for $\hat{d}(\sigma(s),\cdot)$ at $q$ along the reversed segment.
By \eqref{eq:exact-decomposition} this supports $\hat{d}(p,\cdot)$ at $q$.
Completeness of $\gstar$ and properness of the ray imply
$\int_0^\infty u(\sigma)^{-\theta}ds=\infty$.
Thus these supports have arbitrarily small upper Laplacian at $q$,
which proves \eqref{eq:global-comparison} in this case as well.
\end{proof}

The proof uses the curvature condition only on the finite segments to
which the index form is applied. In particular, the approximating
metrics $u_j^{2\beta}g$ are used only to produce limiting curves;
no curvature inequality for $u_j$ is asserted.

\begin{lemma}[Distance from compact boundary]\label{lem:boundary-distance}
Under the curvature, completeness, and boundary assumptions of
Theorem~\ref{thm:boundary}, the function $r=\hat{d}(B,\cdot)$ is continuous
on $M$, positive on $\interior$, zero on $B$, and satisfies
\begin{equation}\label{eq:r-superharmonic}
 \Delta_g r\leq0\quad\text{on }\interior
\end{equation}
in the upper support sense.
\end{lemma}

\begin{proof}
Compactness of $B$ and the local properties of $\hat{d}$ give positivity
and continuity. Suppose first that $r(q)$ is attained by a
$\hat{d}$-minimizing segment $\sigma$ from $p\in B$ to $q$, parametrized by
$\ghat$-arclength $\tau\in[0,L]$, with unit tangent $T$. Since the
initial point is free to move in $B$, the first variation with
tangential initial fields forces $T(0)\perp T_pB$; the segment leaves
$B$ inward, so $T(0)=\widehat\nu(p)$, the inward $\ghat$-unit normal.
In terms of the inward $g$-unit normal this reads
$T(0)=u(p)^{-\beta}\nu(p)$.

\smallskip\noindent
\emph{The upper support with free initial endpoint.}
Let $E_1,\ldots,E_m$ be a $\ghat$-parallel orthonormal frame normal to
$T$; at $\tau=0$ it spans $T_pB$. Let $\xi$ be smooth with
$\xi(L)=1$; the value $\xi(0)$ is now \emph{free}, since the initial
point may move in $B$. Put $V_i=\xi E_i$ and $V_{m+1}=(\tau/L)T$, and
take the auxiliary metric $k$ to be a product metric in a smooth
collar $B\times[0,\delta)$ of $B$, extended smoothly to $M$. Define
$c_z$ and $F$ by \eqref{eq:simultaneous-variation} and
\eqref{eq:support-definition}, as in Lemma~\ref{lem:point-comparison}.
In the product collar the exponential map of $k$ moves the base
point of $\sigma$ within $B$ when the initial vector is tangent to
$B$: indeed $V_i(0)=\xi(0)E_i(0)\in T_pB$ for $i\leq m$ and
$V_{m+1}(0)=0$, so $c_z(0)\in B$ for all small $z$. Every $c_z$ thus
joins $B$ to $\mathcal E(z)$, whence
\[
 F\geq \hat{d}(B,\cdot)=r\quad\text{near }q,\qquad F(q)=L=r(q).
\]
The segment cannot meet $B$ again (discarding the initial subarc
would shorten it), its initial velocity is inward transverse, and
the rest of the segment stays in the interior; so for small $z$ the
curves remain admissible competitors in $M$. As before, $F$ is
$C^2$, and $dF_q=\ghat(T(L),\cdot)$.

\smallskip\noindent
\emph{Second variation with a free initial endpoint.}
The computation leading to \eqref{eq:second-variation-length} applies
verbatim, except that $V_i(0)\neq0$, so the boundary term at
$\tau=0$ survives:
\begin{equation}\label{eq:second-variation-free}
 \delta_j\delta_iF
 =I(V_i^\perp,V_j^\perp)
  +\big\langle\widehat\nabla_{V_j}V_i,\,T\big\rangle_{\ghat}(L)
  -\big\langle\widehat\nabla_{V_j}V_i,\,T\big\rangle_{\ghat}(0),
\end{equation}
where $I$ is the index form in \eqref{eq:hess-equals-index}. The
chain rule \eqref{eq:hess-chain-rule} cancels the terminal term
exactly as before; there is no chain-rule correction at the initial
endpoint, since $F$ depends only on the terminal point through
$\mathcal E$. The initial term therefore persists in
$\Hess_{\ghat}F$ and is computed from the geometry of $B$: for a
tangent field $W$ along $B$,
\[
 \big\langle\widehat\nabla_W W,\,T(0)\big\rangle_{\ghat}
 =\big\langle\widehat\nabla_W W,\,\widehat\nu\big\rangle_{\ghat}
 =-\big\langle W,\,\widehat\nabla_W\widehat\nu\big\rangle_{\ghat}
 =-\widehat{\mathrm{II}}_{\rm in}(W,W),
\]
where $\widehat{\mathrm{II}}_{\rm in}$ is the second fundamental form
of $B$ for the inward normal in $\ghat$, and we used
$\langle W,\widehat\nu\rangle=0$. With $W=V_i(0)=\xi(0)E_i(0)$, the
trace of the initial terms over $i=1,\ldots,m$ contributes
$+\xi(0)^2\widehat H_{\rm in}(p)$ to $\Delta_{\ghat}F(q)$, where
$\widehat H_{\rm in}=\operatorname{tr}\widehat{\mathrm{II}}_{\rm in}$.
As in Lemma~\ref{lem:point-comparison}, substituting
$\Ric_{\ghat}(T,T)\geq-f''-(f')^2/m$ and integrating
$\int_0^L\xi^2 f''\,d\tau$ by parts produces the boundary bracket
$[\xi^2f']_0^L$: the upper endpoint $+\xi(L)^2f'(L)=f'(L)$ cancels
the drift $-f'(L)$ from \eqref{eq:laplace-transform}, while the lower
endpoint now leaves $-\xi(0)^2f'(0)=-\xi(0)^2\,\widehat\nu f(p)$,
since $f'(0)=df(T(0))=\widehat\nu f(p)$ (primes are $d/d\tau$ and
$T(0)=\widehat\nu$). Collecting,
\begin{equation}\label{eq:free-index}
 e^{-2\phi(q)}\Delta_gF(q)
 \leq m\int_0^L\left(\xi'-\frac{f'}m\xi\right)^2d\tau
       +\xi(0)^2\big(\widehat H_{\rm in}
                              -\widehat\nu f\big)(p),
 \qquad \xi(L)=1 ,
\end{equation}
where the square was completed exactly as in
\eqref{eq:index-square}, now keeping the free value $\xi(0)$.

\smallskip\noindent
\emph{The conformal change of the boundary term.}
For tangent $X,Y$ to $B$, the second fundamental forms of $B$ in
$\ghat=e^{2\phi}g$ and in $g$ (same inward unit direction $\nu$,
$\widehat\nu=e^{-\phi}\nu$) satisfy
\[
 \widehat{\mathrm{II}}_{\rm in}(X,Y)
 =e^{\phi}\big(\mathrm{II}_{\rm in}(X,Y)+\nu(\phi)\,g(X,Y)\big).
\]
Taking the trace with respect to $\ghat|_B=e^{2\phi}g|_B$ gives
$\widehat H_{\rm in}=e^{-\phi}\big(H+m\,\nu(\phi)\big)$, and
$\widehat\nu f=e^{-\phi}\nu f$. With $\phi=\beta\log u$ and
$f=m\log u$, the normal derivatives combine to give
\begin{equation}\label{eq:boundary-transform}
 \widehat H_{\rm in}-\widehat\nu f
 =u(p)^{-\beta}\big(H+m\beta\,\nu(\log u)-m\,\nu(\log u)\big)(p)
 =u(p)^{-\beta}\big(H+\beta\,\nu(\log u)\big)(p)
 =u(p)^{-\beta}q_B(p),
\end{equation}
using $m\beta-m=m(\beta-1)=\beta$ and the definition \eqref{eq:boundary-q}
of $q_B$ from Theorem~\ref{thm:boundary}.

\smallskip\noindent
\emph{Choice of cutoff and conclusion.}
In $\ghat$-arclength choose $\xi(\tau)=u(\sigma(\tau))/u(q)$.
Then $\xi(L)=1$, $\xi(0)=u(p)/u(q)$, and $f=m\log u$ gives
$\partial_\tau\xi-(\partial_\tau f)\xi/m=0$, so the integral
in \eqref{eq:free-index} vanishes. Multiplying \eqref{eq:free-index}
by $e^{2\phi(q)}=u(q)^{2\beta}$ and using $2\beta-2=2/(m-1)$,
$2-\beta=\theta$,
\begin{equation}\label{eq:boundary-support}
 \Delta_gF(q)\leq
 u(q)^{2\beta-2}\,u(p)^{2-\beta}\,q_B(p)
 =u(q)^{2/(m-1)}u(p)^\theta q_B(p)\leq0,
\end{equation}
the last step by the hypothesis $q_B\leq0$ on $B$ and $u>0$.
Thus $r$ has a $C^2$ upper support with nonpositive Laplacian at $q$
whenever a curve from $B$ to $q$ realizes $r(q)$.

Compactness of $B$ gives a point $p\in B$ with $\hat{d}(p,q)=r(q)$.
Suppose now that no curve from $B$ to $q$ has length $r(q)$ in $\ghat$.
Use the complete
approximating metrics $\widehat g_j=u_j^{2\beta}g$ and distances
$\hat{d}_j$ from Lemma~\ref{lem:point-comparison}, and write
$r_j(x)=\min_{w\in B}\hat{d}_j(w,x)$. Since $B$ is compact and $\hat{d}_j\to \hat{d}$
locally uniformly on $M\times M$, $r_j\to r$ locally uniformly.
Each $\hat{d}_j$ is proper and geodesic, so there is a $\hat{d}_j$-minimizing
curve $c_j$ from $q$ to $B$, parametrized by $g$-arclength from $q$.

If the $g$-lengths of these curves stay bounded along a subsequence,
reparametrize on $[0,1]$ with constant $g$-speed. Their images lie
in a fixed compact set, so Arzel\`a--Ascoli and lower semicontinuity
of weighted length give a path from $q$ to $B$ of weighted length at
most $\lim r_j(q)=r(q)$, contrary to our assumption. Thus the $g$-lengths
tend to infinity. On every finite parameter interval the curves are
$1$-Lipschitz for $d_g$ and lie in a fixed compact $g$-ball. By
Arzel\`a--Ascoli a diagonal subsequence converges locally uniformly to
a $1$-Lipschitz curve $\sigma:[0,\infty)\to M$ starting at $q$.
For $0\leq s<t$, the subarcs of $c_j$ minimize $\hat{d}_j$ and eventually
lie where $u_j=u$. Local uniform convergence of the distances gives
$$
  \hat{d}(\sigma(s),\sigma(t))
  =\lim_j\hat{d}_j(c_j(s),c_j(t))
  =\int_s^t u(\sigma(v))^\beta\,dv .
$$
The $\ghat$-length of the limit subarc lies between this distance and
the last integral, since $|\dot\sigma|_g\leq1$ almost everywhere.
Thus every subarc minimizes $\hat{d}$, and positivity of $u$ gives
$|\dot\sigma|_g=1$ almost everywhere. Splitting each $c_j$ at $s$
and using $r_j\to r$ gives
\begin{equation}\label{eq:boundary-decomposition}
 r(q)=\hat{d}(q,\sigma(s))+r(\sigma(s))\qquad(s>0).
\end{equation}
If $\sigma(s)\in B$ for some finite $s$, its initial subarc would
realize $r(q)$. Thus $\sigma$ stays in $\interior$, where each subarc
is a $\ghat$-geodesic after reparametrization by $\ghat$-arclength.
It is proper by the
compact return argument of Lemma~\ref{lem:point-comparison}: a
return to a compact set infinitely often would make the weighted
lengths of $[0,t_k]$ bounded while each unit interval contributes at
least $\min_{C_1}u^\beta>0$.

Fix $s>0$ and let $F_s$ be the point-distance upper support for
$\hat{d}(\sigma(s),\cdot)$ at $q$ along the reversed segment
$\sigma|_{[0,s]}$, from Lemma~\ref{lem:point-comparison}. Set
\[
 S_s=F_s+r(\sigma(s)).
\]
By \eqref{eq:boundary-decomposition}, $S_s(q)=r(q)$; and for any $x$,
the triangle inequality gives
$S_s(x)\ge \hat{d}(\sigma(s),x)+r(\sigma(s))\ge r(x)$. Thus $S_s$ is a
$C^2$ upper support for $r$ at $q$, with
\[
 \Delta_gS_s(q)
 \leq\frac{m\,u(q)^{2/(m-1)}}{\displaystyle\int_0^s u(\sigma)^{-\theta}\,dr}.
\]
Completeness of $\gstar$ and properness of $\sigma$ imply
$\int_0^\infty u(\sigma)^{-\theta}\,dr=\infty$, so the right side
tends to $0$ as $s\to\infty$. For every $\varepsilon>0$, choosing
$s$ large gives an upper support with Laplacian at most
$\varepsilon$ at $q$. This is \eqref{eq:r-superharmonic} in the
upper support sense.
\end{proof}

\section{Busemann functions and harmonicity}\label{sec:busemann}

This section constructs the weighted Busemann function that drives the
splitting. We first produce weighted minimizing lines in the two-ended
case and weighted minimizing rays from the boundary in the
compact boundary case. Applying the distance comparisons of the previous
section to the associated Busemann approximations, we then prove that the
limit function $b$ is smooth, harmonic for $\Delta_g$, and satisfies the
eikonal equation $|\nabla b|_g=u^\beta$, with $\nabla b=u^\beta\nu$ on the
boundary in the boundary case. These properties are the input for the
equality analysis of the next section.

\begin{lemma}[Weighted minimizing lines and rays]\label{lem:curves}
Let $(M,g)$ be complete and connected. Every curve $\sigma:\R\to M$
parametrized by $g$-arclength and minimizing $\hat{d}$ on each compact
subarc is proper in both directions. If $M$ is boundaryless with at
least two ends, such a line exists.
If $M$ is noncompact with nonempty compact boundary $B$,
there is a proper ray $\sigma:[0,\infty)\to M$, also parametrized by
$g$-arclength, such that
$\sigma(0)\in B$ and
\begin{equation}\label{eq:free-ray}
 \hat{d}(B,\sigma(s))=L_{\ghat}(\sigma|_{[0,s]})\quad(s\geq0).
\end{equation}
\end{lemma}

\begin{proof}
First consider any $g$-unit-speed curve $\sigma:[0,\infty)\to M$
whose subarcs minimize $\hat{d}$. If $\sigma(t_k)$ belongs to a
compact set $C$ for a sequence $t_k\to\infty$, pass to a subsequence with
$t_{k+1}>t_k+1$. The weighted length of $\sigma|_{[0,t_k]}$ is at most
$\max_{y\in C}\hat{d}(\sigma(0),y)$. For each $i<k$, the interval
$[t_i,t_i+1]$ maps into the compact closed unit $g$-neighborhood $C_1$
of $C$ and contributes at least $\min_{C_1}u^\beta>0$ to this length.
Summing over $i<k$ gives a contradiction. Applying this argument to
both orientations of a weighted line proves the properness assertion.

For existence we adapt the approximation used in the proof of
\cite[Proposition~2.4]{HW26a} to the present coefficient and to the case
with boundary. All curve parameters below are measured in $g$-arclength.
Take the complete metrics $\widehat g_j=u_j^{2\beta}g$
and distances $\hat{d}_j$ constructed in Lemma~\ref{lem:point-comparison}.
The local uniform convergence $\hat{d}_j\to \hat{d}$ holds on $M\times M$,
including boundary points for intrinsic distances.
In particular this convergence allows both endpoints to vary.
\begin{itemize}
    \item

In the boundaryless case choose a compact set $K$ for which $M\setminus K$
has two distinct unbounded components. In each choose $p_j^\pm$ with
$d_g(p_j^\pm,K)\to\infty$, and join them by a $\hat{d}_j$-minimizer $c_j$.
It meets $K$; choose one such point to have parameter zero. Both
parameter lengths from zero to the endpoints tend to infinity, since
each is at least the respective $g$-distance from $K$. For fixed $T$,
the image of $[-T,T]$ lies in the compact closed $g$-neighborhood of
radius $T$ of $K$. The curves are $1$-Lipschitz for $d_g$, so a diagonal
compactness argument gives a locally uniform limit
$\sigma:\R\to M$ which is $1$-Lipschitz for $d_g$.

For fixed $s<t$, the subarcs $c_j|_{[s,t]}$ are themselves minimizing
and are eventually contained where $u_j=u$. Therefore
\begin{align*}
 \hat{d}(\sigma(s),\sigma(t))
 &=\lim_j \hat{d}_j(c_j(s),c_j(t))
 =\lim_j\int_s^t u(c_j(r))^\beta\,dr
 =\int_s^t u(\sigma(r))^\beta\,dr .
\end{align*}
The weighted length of the limit subarc is between this distance and
the last integral: the upper bound follows from $|\dot\sigma|_g\leq1$
almost everywhere. Equality and $u>0$ thus give $|\dot\sigma|_g=1$
almost everywhere and show that every limit subarc minimizes $\hat{d}$.
Locally in the interior, reparametrization by its positive weighted
speed makes it a $\ghat$-geodesic.

\item In the boundary case take $p_j$ with $d_g(B,p_j)\to\infty$ and a
$\hat{d}_j$-minimizing curve from the compact set $B$ to $p_j$, starting at
$c_j(0)\in B$. The same argument, with $K=B$ and $[0,T]$ in place
of $[-T,T]$, gives a $g$-unit-speed limit on $[0,\infty)$.
Write $r_j(x)=\min_{w\in B}\hat{d}_j(w,x)$ and $r(x)=\hat{d}(B,x)$. Compactness
of $B$ and uniform convergence on $B\times C$ for every compact $C$
give $r_j\to r$ locally uniformly. Every initial subarc realizes the
distance from $B$ (otherwise one could shorten the full curve), hence
$$
 r(\sigma(s))=\lim_j r_j(c_j(s))
 =\int_0^s u(\sigma(r))^\beta\,dr
 =L_{\ghat}(\sigma|_{[0,s]}).
$$
Thus $\sigma(0)\in B$, and the right side is positive for $s>0$:
the limiting ray never meets $B$ again. Its short initial subarc is
normal to $B$: it minimizes distance from $B$, and in a tubular
neighborhood the normal collar and Gauss's lemma give precisely the
normal segments as minimizers. Conformal change preserves
orthogonality, so its initial $g$-unit tangent is the inward normal.
\end{itemize}
Every subarc of the boundary ray also minimizes $\hat{d}$, since
shortening it would contradict \eqref{eq:free-ray}. The initial
properness argument therefore applies to this ray as well.
\end{proof}

\begin{proposition}[Harmonic functions with prescribed gradient]
\label{prop:harmonic}
Under the hypotheses of Theorem~\ref{thm:closed}, there is a smooth function $b$ on $M$
such that
\begin{equation}\label{eq:harmonic-eikonal}
 \Delta_gb=0,\qquad |\nabla b|_g=u^\beta>0.
\end{equation}
On the prescribed line $\sigma$, this function satisfies $b(\sigma(0))=0$
and $\dot\sigma=\nabla b/|\nabla b|_g$.
Under the hypotheses of Theorem~\ref{thm:boundary}, $b=r=\hat{d}(B,\cdot)$ is smooth up to $B$,
satisfies \eqref{eq:harmonic-eikonal} on $\interior$, and
$\nabla b=u^\beta\nu$ on $B$.
\end{proposition}

\begin{proof}
\emph{The boundaryless case.}
Let $\sigma$ be the line given in Theorem~\ref{thm:closed}, let $o=\sigma(0)$,
and put $p_s^\pm=\sigma(\pm s)$ and
$b_s^\pm(x)=\hat{d}(o,p_s^\pm)-\hat{d}(x,p_s^\pm)$ for $s\geq0$.
For $t\geq s$, the minimizing property gives
$\hat{d}(o,p_t^\pm)=\hat{d}(o,p_s^\pm)+\hat{d}(p_s^\pm,p_t^\pm)$.
The triangle inequality for $\hat{d}(x,p_t^\pm)$ therefore gives
$b_t^\pm(x)\geq b_s^\pm(x)$. Also
$$
 |b_s^\pm(x)|\leq \hat{d}(o,x),\qquad
 |b_s^\pm(x)-b_s^\pm(y)|\leq \hat{d}(x,y).
$$
Thus the limits $b^\pm$ are finite, are $1$-Lipschitz for $\hat{d}$, and
convergence is uniform on compact sets (apply the Lipschitz bound
on a finite $\hat{d}$-net). No limit of the individual distances
$\hat{d}(o,p_s^\pm)$ is required to be infinite.

For a compact set $K$, and $s$ large enough that $p_s^\pm\notin K$,
Lemma~\ref{lem:point-comparison} gives the lower-support bound
$$
 \Delta_g b_s^\pm(x)\geq
 -\frac{m\,u(x)^{2/(m-1)}}{d_{\gstar}(x,p_s^\pm)}
 \qquad(x\in K).
$$
The complete length metric of $\gstar$ is proper and has the original
manifold topology. By Lemma~\ref{lem:curves} both tails of $\sigma$ are proper, so
$d_{\gstar}(o,p_s^\pm)\to\infty$. The triangle inequality bounds the
denominator below by $d_{\gstar}(o,p_s^\pm)-\max_K d_{\gstar}(o,\cdot)$.
Consequently the error tends uniformly to zero on $K$.
Passing to the limit gives $\Delta_g b^\pm\geq0$ in distributions;
equivalently one may integrate the preceding bound against each
nonnegative compactly supported test function before passing to the limit.
For each $s$,
\begin{align*}
 b_s^+(x)+b_s^-(x)
 &=\hat{d}(p_s^+,p_s^-)-\hat{d}(x,p_s^+)-\hat{d}(x,p_s^-)\leq0 ,
\end{align*}
and both functions vanish at $o$. The continuous subharmonic function
$b^++b^-$ attains its maximum there. Connectedness and the strong maximum
principle give $b^++b^-=0$ everywhere. Their distributional Laplacians
are nonnegative and sum to zero, hence each vanishes. Elliptic
regularity makes $b=b^+$ smooth.

\smallskip\noindent
\emph{The boundary case: keeping the comparison in the interior.}
Write $p_s=\sigma(s)$ for the free ray, $r=\hat{d}(B,\cdot)$, and set
$b_s(x)=r(p_s)-\hat{d}(x,p_s)$. For $0\leq s\leq t$,
\eqref{eq:free-ray} gives $r(p_t)=r(p_s)+\hat{d}(p_s,p_t)$.
Thus $b_s$ again increases and is $1$-Lipschitz for $\hat{d}$.
Since $-\hat{d}(x,p_0)=b_0(x)\leq b_s(x)\leq r(x)$, it converges locally
uniformly to a continuous function $b\leq r$. On the ray,
$b_t(p_s)=r(p_s)$ for $t\geq s$; in particular
$b(p_s)=r(p_s)>0$ for $s>0$.

Let $\Omega=\{b>0\}\subset\interior$ and take $U\Subset\Omega$ open.
Choose $\delta>0$ with $b>3\delta$ on $\overline U$.
For all sufficiently large $s$, $b_s>2\delta$ on $\overline U$
and $p_s\notin\overline U$. \emph{Fix such an $s$ first.}
The locally uniform convergence $\hat{d}_j\to \hat{d}$ from
Lemma~\ref{lem:point-comparison} implies, for all sufficiently large $j$,
$$
 \hat{d}_j(x,p_s)<r(p_s)-\delta\qquad(x\in\overline U).
$$
If $z$ lies on a $\hat{d}_j$-minimizing curve from $x$ to $p_s$, then
$$
 r(p_s)-r(z)\leq \hat{d}(z,p_s)\leq \hat{d}_j(z,p_s)\leq \hat{d}_j(x,p_s).
$$
Hence $r(z)>\delta$ along the entire curve. In particular all of its
locally convergent pieces have limits in $\{r\geq\delta\}\subset\interior$.

We spell out the point-distance comparison at each $x\in U$,
since the global statement of Lemma~\ref{lem:point-comparison}
was made only for manifolds without boundary. If the approximating
curves have a subsequence of bounded $g$-length, the compactness
argument of that lemma gives an interior $\hat{d}$-minimizer from $x$ to
$p_s$. The reversed segment in \eqref{eq:segment-comparison}
has $\gstar$-length at least $d_{\gstar}(x,p_s)$.
Otherwise one obtains a proper interior limiting ray $\eta$ from $x$.
Splitting the approximating distances at a fixed parameter $T>0$ gives
$\hat{d}(x,p_s)=\hat{d}(x,\eta(T))+\hat{d}(\eta(T),p_s)$.
An upper support for $\hat{d}(\eta(T),\cdot)$ at $x$, plus the constant
$\hat{d}(\eta(T),p_s)$, is then an upper support for $\hat{d}(p_s,\cdot)$.
Its Laplacian can be bounded using the finite reversed segment:
the denominator in \eqref{eq:segment-comparison} tends to infinity
as $T\to\infty$, since $\eta$ is proper and $\gstar$ is complete.
Both alternatives therefore give, in the lower-support sense,
$$
 \Delta_g b_s(x)\geq
 -\frac{m\,u(x)^{2/(m-1)}}{d_{\gstar}(x,p_s)}\qquad (x\in U).
$$
Curvature was used only on interior segments for the original metric.
Now let $s\to\infty$. Properness of the free ray gives uniform
divergence of the denominator on $\overline U$ as above.
It follows that $\Delta_g b\geq0$ on $\Omega$.

On $\Omega$, $b-r\leq0$ is subharmonic by
Lemma~\ref{lem:boundary-distance}. It is zero on the positive part
of the ray, so the strong maximum principle gives $b=r$ on the
connected component $U_0$ of $\Omega$ containing that part.
If $x\in\overline{U_0}\cap\interior$, then $b(x)=r(x)>0$ by
continuity and positivity of $r$ in the interior. A connected
neighborhood of $x$ in $\Omega$ intersects $U_0$, so $x\in U_0$.
Thus $U_0$ is both open and closed in $\interior$, which is connected
(a path in $M$ between interior points can be pushed off $B$ in a collar).
We conclude that $b=r$ throughout the interior, and hence on the
boundary by continuity. It is both subharmonic and superharmonic in
the interior, and thus smooth and harmonic there.

In a smooth boundary chart $r$ is locally Lipschitz up to $B$:
it is $1$-Lipschitz for $\hat{d}$ and the weight is bounded on a smaller
coordinate half-ball. Its trace on $B$ is zero. On such a half-ball
the weak harmonic equation extends by density from interior test
functions to $H^1$ test functions of zero trace with support away from
the artificial boundary. Local Dirichlet boundary regularity for the
smooth metric $g$ and the zero boundary data now gives smoothness up to $B$
\cite[Theorem~6.19]{GilbargTrudinger}.

\smallskip\noindent
\emph{The gradient identity in both cases.}
Use $p_s=p_s^+$ and $b_s=b_s^+$ in the boundaryless case.
Fix any interior point $x$. As $b$ is $1$-Lipschitz for $\hat{d}$, comparison
along short smooth curves through $x$ gives $|\nabla_{\ghat}b(x)|\leq1$.
For the opposite inequality choose a $\ghat$-normal ball $V$ centered
at $x$ of radius $\ell>0$ with compact closure in a larger interior
normal neighborhood. Gauss's lemma gives
$\hat{d}_V(x,y)=\ell$ for $y\in\partial V$, where $\hat{d}_V$ measures length in
$\overline V$. For large $s$, $p_s\notin\overline V$, and
\begin{equation}\label{eq:dynamic}
 \hat{d}(x,p_s)=\ell+\min_{y\in\partial V}\hat{d}(y,p_s).
\end{equation}
Indeed, every path to $p_s$ first meets $\partial V$ after length
at least $\ell$. Conversely, follow a radial segment to a minimizer
of $\hat{d}(\cdot,p_s)$ on the compact sphere and then follow an
arbitrarily close competitor for the remaining distance.
In both normalizations of $b_s$, \eqref{eq:dynamic} gives
$\max_{\partial V} b_s=b_s(x)+\ell$. Uniform convergence on
$\partial V\cup\{x\}$ yields $y\in\partial V$ with
$b(y)-b(x)=\ell$. Along the radial $\ghat$-unit-speed segment
$c:[0,\ell]\to\overline V$ from $x$ to $y$ the Lipschitz bounds give
$$
 b(y)-(\ell-t)\ \leq\ b(c(t))\ \leq\ b(x)+t.
$$
The two bounds coincide. Differentiating at zero gives
$db_x(\dot c(0))=1$, whence $|\nabla_{\ghat}b|=1$. By conformal
change $|\nabla_g b|=u^\beta$ at every interior point.
In the boundaryless case, for $t\in\R$ and $s>|t|$ the minimizing
property of the oriented line gives
$$
 b_s^+(\sigma(t))=\hat{d}(\sigma(0),\sigma(s))-
 \hat{d}(\sigma(t),\sigma(s))=\int_0^t u(\sigma(r))^\beta\,dr .
$$
Passing to the limit gives $b(\sigma(0))=0$ and the same identity for
$b(\sigma(t))$. Differentiation gives
$g(\nabla b,\dot\sigma)=u^\beta$ almost everywhere on the line.
Since $|\nabla b|_g=u^\beta$ and $|\dot\sigma|_g=1$, equality in
Cauchy--Schwarz gives $\dot\sigma=\nabla b/|\nabla b|_g$ almost
everywhere. Thus $\sigma$ is an integral curve of this smooth unit
field, and the identity holds for every $t$.

In the boundary case the identity $|\nabla_g b|=u^\beta$ extends to $B$ by continuity.
Since $b=0$ on $B$, all its tangential derivatives vanish there;
since $b\geq0$ on the inward side, $\nu b\geq0$. The nonzero
gradient norm then gives $\nabla_g b=u^\beta\nu$ on $B$.
\end{proof}

\begin{remark}[Relation to the classical Busemann function]
\label{rem:busemann}
The function $b$ of Proposition~\ref{prop:harmonic} is a weighted Busemann
function. In the Cheeger--Gromoll proof \cite{CheegerGromoll71}, a line
$\sigma$ in a manifold of nonnegative Ricci curvature yields the Busemann
functions $b^\pm(x)=\lim_{s\to\infty}\big(d(o,\sigma(\pm s))-d(x,\sigma(\pm s))\big)$;
each is subharmonic by the Laplacian comparison, $b^++b^-\leq0$ with
equality at $o$, and the maximum principle forces $b^++b^-=0$, so
$b=b^+$ is harmonic with $|\nabla b|=1$ and its gradient flow splits off a
line. Here the Riemannian distance $d$ is replaced by the weighted
distance $\hat{d}$ of the conformal metric $\ghat=u^{2\beta}g$.
Lemma~\ref{lem:point-comparison} bounds $\Delta_g b_s^\pm$ from below
by a negative error tending to zero uniformly on compact sets.
Passing to the limit gives $\Delta_g b^\pm\geq0$ in distributions.
The rest of the argument is parallel: $b^++b^-\leq0$ with equality at $o$
again forces harmonicity, but the eikonal identity becomes
$|\nabla b|_g=u^\beta$ instead of $|\nabla b|=1$. The Bochner--Kato
equality of Lemma~\ref{lem:equality} then yields a warped product rather
than a Riemannian product; when $u$ is constant this reduces to the
classical splitting. The same weighted Busemann construction appears in
\cite{HW26a}. In the compact boundary case, $q_B\leq0$ gives
$\Delta_g r\leq0$ for $r=\hat{d}(B,\cdot)$. On $\{b>0\}$ the function
$b-r$ is subharmonic, is nonpositive, and vanishes along the interior
of the free ray. The strong maximum principle gives $b=r$ on the
component containing the ray. By continuity and $r>0$ in $\interior$,
this component is also closed in the connected interior, so $b=r$
throughout $M$.
\end{remark}

\section{Equality, global flow, and curvature}\label{sec:rigidity}

This section turns the comparison into a rigid splitting. We analyze the
equality case of the Bochner--Kato inequality to show that the level sets
of $b$ are totally umbilical and that the integral curves of their unit
normal $N=\nabla b/|\nabla b|$ are unit-speed geodesics.
A global flow argument then writes the metric as a warped
product $dt^2+a(t)^2h$, identifies the warping function with $u$, and
derives the fiber Ricci bound. Combining these ingredients with the
boundary analysis proves Theorems~\ref{thm:closed} and~\ref{thm:boundary}.

\begin{lemma}[The equality mechanism]\label{lem:equality}
Suppose \eqref{eq:main-curvature} and \eqref{eq:harmonic-eikonal}
hold. Put $v=|\nabla b|=u^\beta$ and $N=\nabla b/v$.
Then
\begin{equation}\label{eq:normal-geometry}
 \nabla_NN=0,\qquad
 X(u)=0,\qquad
 \nabla_XN=-\frac{N(\log u)}{n-2}\,X
 \quad(X\perp N),
\end{equation}
and
$\Ric(N,N)=\rho=\beta\Delta u/u$.
\end{lemma}

\begin{proof}
Since $v=u^\beta$, $\nabla v=\beta u^{\beta-1}\nabla u$ and
$\Delta v=\beta u^{\beta-1}\Delta u+\beta(\beta-1)u^{\beta-2}|\nabla u|^2$.

Using $(\beta-1)/\beta=1/m$ and $|\nabla v|^2=\beta^2 u^{2\beta-2}|\nabla u|^2$,
\[
 v\Delta v=\beta\frac{\Delta u}{u}v^2+\beta(\beta-1)u^{2\beta-2}|\nabla u|^2
          =\beta\frac{\Delta u}{u}v^2+\frac1m|\nabla v|^2.
\]
The Bochner formula for the harmonic function $b$ reads
\[\tfrac12\Delta v^2=|\Hess b|^2+\Ric(\nabla b,\nabla b)
   =|\Hess b|^2+v^2\Ric(N,N).\]
Since $\tfrac12\Delta v^2=v\Delta v+|\nabla v|^2$,
substituting the expression for $v\Delta v$ and rearranging gives
\begin{equation}\label{eq:bochner-equality}
 0=\left(|\Hess b|^2-\frac{m+1}{m}|\nabla v|^2\right)
    +v^2\left(\Ric(N,N)-\beta\frac{\Delta u}{u}\right).
\end{equation}
The second summand is $\geq0$ by \eqref{eq:main-curvature}. For the first,
take an orthonormal frame $N,e_1,\ldots,e_m$ and write
\[
 h_{00}=\Hess b(N,N),\qquad
 z_i=\Hess b(N,e_i),\qquad
 \mathcal B_{ij}=\Hess b(e_i,e_j).
\]
Then $|\Hess b|^2=h_{00}^2+2|z|^2+|\mathcal B|^2$. Harmonicity gives
$\operatorname{tr}\mathcal B=-h_{00}$, and $e_i(v)=\Hess b(N,e_i)=z_i$
yields $|\nabla v|^2=h_{00}^2+|z|^2$. Hence
\[
 |\Hess b|^2-\frac{m+1}{m}|\nabla v|^2
 =|\mathcal B|^2-\frac1m h_{00}^2+\frac{m-1}{m}|z|^2 .
\]
Since $\operatorname{tr}\mathcal B=-h_{00}$, completing the square gives
\[
 \left|\mathcal B+\frac{h_{00}}m I\right|^2
 =|\mathcal B|^2+\frac{2h_{00}}m\operatorname{tr}\mathcal B+\frac{h_{00}^2}{m}
 =|\mathcal B|^2-\frac{2h_{00}^2}{m}+\frac{h_{00}^2}{m}
 =|\mathcal B|^2-\frac1m h_{00}^2 .
\]
Thus
\begin{equation}\label{eq:kato-defect}
 |\Hess b|^2-\frac{m+1}{m}|\nabla v|^2
 =\left|\mathcal B+\frac{h_{00}}m I\right|^2+\frac{m-1}{m}|z|^2\geq0 .
\end{equation}
Both summands in \eqref{eq:bochner-equality} are therefore nonnegative and
sum to zero, so each vanishes. From \eqref{eq:kato-defect},
$z=0$ and $\mathcal B=-(h_{00}/m)I$.
Since $e_i(v)=z_i$, the function $u$ is constant along every connected
level hypersurface. Moreover
$g(\nabla_NN,e_i)=z_i/v=0$ and, for tangent $X$,
\[
 \nabla_XN=-\frac{N(v)}{mv}X
          =-\frac{N(\log u)}{m-1}X.
\]
The second term of \eqref{eq:bochner-equality} vanishes, and
$\Ric(N,N)\geq\rho\geq\beta\Delta u/u$ gives all claimed equalities.
\end{proof}

\begin{proof}[Proof of Theorem~\ref{thm:closed}]
Use Proposition~\ref{prop:harmonic} and Lemma~\ref{lem:equality}.
Write $x_0=\sigma(0)$ and let $\Sigma$ be the connected component of
$b^{-1}(0)$ containing $x_0$. It is an embedded hypersurface without
boundary because $db$ never vanishes.
As a connected component of the closed level set, $\Sigma$ is closed in $M$.
An intrinsic Cauchy sequence converges in $M$ to a point of $\Sigma$,
and a submanifold chart there gives convergence for the induced metric
$h$. Thus $h$ is complete. The unit field $N$ has a complete flow
$\Phi_t$ by completeness of $g$. The tangential derivatives of $v$
vanish, so $dv=N(v)N^\flat$ and
$dN^\flat=-v^{-2}dv\wedge db=0$. Since $N^\flat(N)=1$, Cartan's
formula gives $\mathcal L_NN^\flat=0$.
Define $F:\R\times\Sigma\to M$ by $F(t,x)=\Phi_t(x)$.
Since $\dot\sigma=N$ by Proposition~\ref{prop:harmonic}, uniqueness
of integral curves gives $F(t,x_0)=\sigma(t)$ for every $t\in\R$.
Then $F^*N^\flat=dt$: the $t$-derivative maps to $N$ and
$N^\flat$ vanishes on vectors transported from $T\Sigma$. These
transported vectors remain linearly independent since $\Phi_t$ is
a diffeomorphism. Together with $N$ they give full rank.
For each $t$, the connected slice $F(\{t\}\times\Sigma)$ is tangent
to $N^\perp$, so $b$ and $u$ are constant on that slice.
Write these values as $\mathfrak b(t)$ and $U(t)$. Then
$\mathfrak b'(t)=U(t)^\beta>0$.
For transported tangent fields $J_X,J_Y$ we have $[N,J_X]=[N,J_Y]=0$.
Metric compatibility and \eqref{eq:normal-geometry} therefore give
$$
 \partial_t\,g(J_X,J_Y)
 =g(\nabla_{J_X}N,J_Y)+g(J_X,\nabla_{J_Y}N)
 =-\frac{2U'(t)}{(m-1)U(t)}g(J_X,J_Y).
$$
The coefficient is constant on each slice. Integration from $t=0$
therefore gives
\begin{equation}\label{eq:flow-pullback}
 F^*g=dt^2+a(t)^2h,\qquad
 a(0)=1,\qquad
 \frac{a'}a=-\frac1{m-1}\frac{U'}U .
\end{equation}
In particular $U=U(0)a^{-(m-1)}$.

The metric in \eqref{eq:flow-pullback} is complete: on each bounded
$t$-interval, $a$ has a positive lower bound, and $h$ is complete;
a curve with unbounded $t$ has infinite length.
To see surjectivity, join $F(0,x_0)$ to a given $q\in M$ by a
$g$-geodesic segment, using Hopf--Rinow. Lift its initial velocity
by $dF_{(0,x_0)}^{-1}$. Completeness of the source provides the
geodesic with this initial data for the full time interval. Its image
has the same geodesic initial data in $M$, hence ends at $q$.
If $F(t,x)=F(s,y)$, strict monotonicity of $\mathfrak b(t)$ gives $t=s$, and
uniqueness of the flow gives $x=y$.
Thus $F$ is a global isometry on $M$ itself.
The curvature assertions follow from Proposition~\ref{prop:converse} below.
\end{proof}

\begin{proof}[Proof of Theorem~\ref{thm:boundary}]
Now $b=r$ and $N|_B=\nu$ by Proposition~\ref{prop:harmonic}.
Transversality and the inverse function theorem give an initial flow
collar (of uniform width since $B$ is compact).
The equations in Lemma~\ref{lem:equality} extend continuously to $B$,
so $u|_B=C$ by connectedness. The forward flow from $B$ exists for
all $t\geq0$: its speed is one, it cannot escape in finite time
by completeness, and it cannot return to the boundary because
$b$ strictly increases along its trajectories.
Let $F(t,p)=\Phi_t(p)$ denote this flow. On $B$, $N^\flat$ vanishes
on tangent vectors, and $\mathcal L_NN^\flat=0$ as in the boundaryless
case. Thus $F^*N^\flat=dt$. The slices are connected and perpendicular
to $N$, so $b\circ F=\mathfrak b(t)$ and $u\circ F=U(t)$ with $\mathfrak b'(t)=U(t)^\beta>0$.
For tangent vectors transported by the flow, Lemma~\ref{lem:equality}
gives $\partial_t(F_t^*g)=-2(U'/U)(F_t^*g)/(m-1)$, where
$F_t(p)=F(t,p)$. The initial value is $h=g|_B$; integration gives
$F^*g=dt^2+a^2h$ and $U=Ca^{-(m-1)}$, with $a(0)=1$.
This makes $F$ a local isometry. If two images coincide, strict
monotonicity of $\mathfrak b(t)$ gives the same time, and uniqueness for
the flow then gives the same point of $B$.

For surjectivity, join a point in the initial collar to any prescribed
interior point by a finite-length piecewise smooth path in $\interior$.
Lift this path through the local isometry. Along a lift,
$|dt|\leq ds$, so the $t$ coordinate stays bounded.
The other coordinate lies in the compact $B$. The length of every
tail of the lift equals that of the corresponding tail of the original
path, so the lift is Cauchy at any finite terminal time.
It therefore has a limiting point in a compact $t$-slab and extends there.
It cannot terminate at $t=0$, because its projected point is in
$\interior$ whereas $F(\{0\}\times B)=B$.
The entire path lifts, proving surjectivity. The normal collar already
covers $B$, so the isometry is global including the boundary.
Finally $H=ma'/a$ and $\nu\log u=-(m-1)a'/a$ at $t=0$ give
$q_B=0$. The remaining statements follow below.
\end{proof}

\begin{proof}[Proof of Corollary~\ref{cor:two-ends}]
Since $M$ has at least two ends, Lemma~\ref{lem:curves} gives a weighted
line for the specified $u$. Apply Theorem~\ref{thm:closed} to this line.
If the connected fiber $\Sigma$ were noncompact, $\R\times\Sigma$ would
have only one end. Indeed, take a compact exhaustion $K_j$ of $\Sigma$.
The upper and lower slabs of the complement of $[-j,j]\times K_j$
are each connected and are joined by a vertical line through any
$x\in\Sigma\setminus K_j$. Every remaining point of the complement
can be joined vertically to one of those slabs. This contradicts the
end hypothesis. Hence $\Sigma$ is compact, and the complements of
$[-j,j]\times\Sigma$ have exactly two unbounded components.
\end{proof}

\noindent\emph{The supersolutions in Corollary~\ref{cor:three}.}
Fix the splitting with compact fiber given by Corollary~\ref{cor:two-ends} for a positive solution
$u_0=C_0/a$ in dimension three. Let $w$ be any positive $C^2$ supersolution
of $L_2$ and set $h_1=w/u_0$. Since $L_2u_0=0$, direct expansion gives
$\operatorname{div}(u_0^2\nabla h_1)=-(u_0/2)L_2 w\leq0$.
Testing with $\eta^2/h_1$ and applying Cauchy--Schwarz yields
$$
 \int_M\eta^2u_0^2|\nabla\log h_1|^2
 \leq 4\int_M u_0^2|\nabla\eta|^2
 \qquad(\eta\in C_c^\infty(M)).
$$
Take $\eta(t)=\chi(t/R)$ with $\chi\in C_c^\infty(\R)$ equal to one
on $[-1,1]$. The fiber is compact and $u_0^2dV_g=C_0^2\,dt\,dV_h$,
so the right-hand side is
$4C_0^2\vol_h(\Sigma)R^{-1}\int_\R|\chi'|^2\to0$.
On each fixed compact set $\eta=1$ for all large $R$, whence
$\nabla\log h_1=0$. Connectedness now gives $w=C/a$ in the
\emph{same} splitting coordinates.

\begin{proposition}[Warped curvature and converse]\label{prop:converse}
Let $h$ be a metric on an $m$-manifold, let $a>0$, and set
$g=dt^2+a^2h$ and $u=Ca^{-(m-1)}$.
Then \eqref{eq:main-curvature} holds if and only if
\begin{equation}\label{eq:converse-fiber}
 \Ric_h\geq(m-1)\big((a')^2-aa''\big)h
\end{equation}
for every $t$ in the base interval. In that case
$\rho=-ma''/a=\beta\Delta u/u$.
\end{proposition}

\begin{proof}
For vectors $X,Y$ lifted from the fiber, the warped curvature formulas are
\begin{align}
 \Ric_g(\partial_t,\partial_t)&=-m\frac{a''}{a},&
 \Ric_g(\partial_t,X)&=0,\label{eq:warped-ricci}\\
 \Ric_g|_{T\Sigma}
 &=\Ric_h-\big(aa''+(m-1)(a')^2\big)h.\nonumber
\end{align}
They follow from
$\nabla_X\partial_t=(a'/a)X$ and
$\nabla_XY=\nabla^h_XY-aa'h(X,Y)\partial_t$
for $t$-independent fiber fields.
Write $k=a'/a$. From $u=Ca^{-(m-1)}$,
\[
 \frac{u'}u=-(m-1)\frac{a'}a=-(m-1)k,\qquad
 \frac{u''}u=\Big(\frac{u'}u\Big)'+\Big(\frac{u'}u\Big)^2
          =-(m-1)k'+(m-1)^2k^2 .
\]
The Laplacian of a radial function is $\Delta u=u''+mk\,u'$, so
\begin{equation}\label{eq:u-laplace-warped}
 \frac{\Delta u}{u}
 =\frac{u''}u+mk\frac{u'}u
 =-(m-1)k'+(m-1)^2k^2-m(m-1)k^2
 =-(m-1)(k'+k^2)
 =-(m-1)\frac{a''}a ,
\end{equation}
since $k'+k^2=a''/a$. Hence
$\beta\Delta u/u=-\beta(m-1)a''/a=-ma''/a=\Ric_g(\partial_t,\partial_t)$,
using $\beta(m-1)=m$.

The tangential Ricci tensor is
$\Ric_h-\big(aa''+(m-1)(a')^2\big)h$, and the radial eigenvalue is
$-ma''/a$. Since the induced metric on each slice is $a^2h$, all
tangential Ricci eigenvalues are at least the radial one precisely when
\[
 \Ric_h-\big(aa''+(m-1)(a')^2\big)h\geq -m\frac{a''}{a}\,a^2 h
        =-ma''a\,h .
\]
Rearranging,
\[
 \Ric_h\geq\big(aa''+(m-1)(a')^2-maa''\big)h
        =(m-1)\big((a')^2-aa''\big)h ,
\]
which is \eqref{eq:converse-fiber}. In that case
$\rho=-ma''/a=\beta\Delta u/u$. This proves the proposition and completes
the curvature parts of both main theorems.
\end{proof}

\section{Nonconstant warped examples}\label{sec:models}

This section shows that the conclusions of the main theorems allow
nonconstant warping. The first model satisfies every hypothesis in all
dimensions $n\geq3$, and its restriction to a half line gives the
compact boundary case. The hyperbolic examples also exhibit a one-ended
manifold with a weighted minimizing line of finite weighted length.

\begin{example}[A model with complete $\gstar$ in every dimension]
\label{ex:nontrivial}
Let $h$ be the unit round metric on $S^{n-1}$ and set
\begin{equation}\label{eq:nontrivial-model}
 a(t)=\sqrt{1+t^2},\qquad
 g=dt^2+(1+t^2)h,\qquad
 u(t)=(1+t^2)^{-(n-2)/2}.
\end{equation}
On $\R\times S^{n-1}$ this is complete and two-ended.
Since $u\leq1$, $\gstar$ is complete.
From $a=(1+t^2)^{1/2}$,
\[
 a'=\frac{t}{(1+t^2)^{1/2}},\qquad
 a''=\frac{1}{(1+t^2)^{3/2}},\qquad
 (a')^2-aa''=\frac{t^2}{1+t^2}-\frac{1}{1+t^2}=\frac{t^2-1}{1+t^2}\leq1 .
\]
By Proposition~\ref{prop:converse}, $\rho=-ma''/a=-(n-1)/(1+t^2)^2$ and
$L_\beta u=0$:
\begin{equation}\label{eq:model-calculation}
 (a')^2-aa''=\frac{t^2-1}{1+t^2}\leq1,\qquad
 \rho=-\frac{n-1}{(1+t^2)^2},\qquad L_\beta u=0 .
\end{equation}
The fiber inequality \eqref{eq:converse-fiber} follows from
$\Ric_h=(n-2)h$ and $(a')^2-aa''\leq1$.
Restricting to $t\geq0$ gives the compact boundary example, with
$a(0)=1$, $H=0$, and $\nu u=0$.
This example has no alternative splitting with a Riemannian line
factor: every plane containing $\partial_t$ has curvature
$-a''/a=-(1+t^2)^{-2}<0$. Any candidate line direction lies in
such a plane, whereas all planes containing a product-line direction
have zero sectional curvature.
\end{example}

\begin{example}[Hyperbolic models]\label{ex:hyperbolic}
Let $(\Sigma,h)$ be a closed hyperbolic surface and take
$a(t)=\cosh t$, $u(t)=\sech t$ on $\R\times\Sigma$.
Then $g=dt^2+\cosh^2t\,h$ has sectional curvature $-1$,
$\rho=-2$, and $(-2\Delta+\rho)u=0$.
No line factor is possible since all sectional curvatures are negative.
Equation~\eqref{eq:ground-state} gives $Q_2\geq0$, and the spectral bottom is zero: $u^2dV_g=dt\,dV_h$, and testing
with $u\chi(t/R)$ in \eqref{eq:ground-state} gives Rayleigh quotients
tending to zero.

The same warping on a noncompact fiber illustrates the line hypothesis.
For $n\geq3$, write hyperbolic space in normal coordinates about a
totally geodesic hyperplane as
$$
 \mathbb H^n=\R\times\mathbb H^{n-1},\qquad
 g=dt^2+\cosh^2t\,h_{\mathbb H^{n-1}},\qquad u=(\cosh t)^{-(n-2)} .
$$
Here $\rho=\beta\Delta u/u=-(n-1)$ by \eqref{eq:u-laplace-warped},
and $\gstar=(\cosh t)^{2(n-3)}g\geq g$ is complete. The function
$F(t)=\int_0^t(\cosh s)^{-(n-1)}\,ds$ satisfies $|dF|_g=u^\beta$.
Hence every curve has weighted length at least the absolute change of
$F$, with equality on each radial segment. Every radial line is thus a
weighted minimizing line of finite total $\ghat$-length, although
$\mathbb H^n$ has only one end.
\end{example}

\section{Counterexamples and the role of conformal completeness}
\label{sec:counterexamples}

In this section, we construct counterexamples to local warped product
splitting over an interval under spectral Ricci bounds in dimensions
$n\geq4$. At $\beta=(n-1)/(n-2)$, a torus example shows that the
completeness of $g_*$ cannot be omitted in either the two-ended or the compact boundary setting. At $c=4/(n-1)$, we construct a metric on
$\mathbb R\times S^{n-1}$ for which local splitting fails even though $\lambda_1(L_c)=0$. A compactly supported perturbation of the potential on this manifold also yields a subcritical Schr\"odinger operator with zero spectral bottom.

Both boundaryless models are complete and two-ended since bounded closed $t$-slabs
are compact and curve length dominates the total variation of $t$.
To rule out local splitting, we use the fact that in $dr^2+a(r)^2h$
the unit radial vector $T$ is a Ricci eigenvector and
$J_T=R(\cdot,T)T=-(a''/a)\operatorname{Id}$ on $T^\perp$.

\subsection{Failure without conformal completeness at \texorpdfstring{$\beta$}{beta}}

\begin{proposition}[An exact counterexample at $\beta$]
\label{prop:beta-counterexample}
For every $n\geq4$, there is a complete two-ended manifold $(M^n,g)$
and a smooth $u>0$ satisfying
\begin{equation}\label{eq:beta-counterexample-zero}
 L_\beta u=0,\qquad \lambda_1(L_\beta)=0,
 \qquad \beta=\frac{n-1}{n-2},
\end{equation}
such that $g_*$ is incomplete and no neighborhood admits a warped-product
representation with one-dimensional base. The same obstruction occurs
on a complete manifold with compact connected boundary satisfying
$H+\beta\nu(\log u)<0$.
\end{proposition}

\begin{proof}
Put $m=n-1>2$ and take
\begin{gather}\label{eq:counterexample-torus-parameters}
 \alpha=\frac{m}{4(m-1)},\qquad b=\frac{1-\alpha}{m-1},
 \qquad a_1=\alpha,\quad a_2=\cdots=a_m=b,\\
 M=\R\times\mathbb T^m,\qquad
 g=dt^2+\sum_i e^{2a_it}d\vartheta_i^2,\qquad u=e^{-t/2}.
 \label{eq:counterexample-torus-metric}
\end{gather}
We have $\sum_i a_i=1$ and
$\alpha-b=(m-2)^2/[4(m-1)^2]>0$.
In the frame $E_0=\partial_t$, $E_i=e^{-a_it}\partial_{\vartheta_i}$,
the identities $\nabla_{E_0}E_i=0$,
$\nabla_{E_i}E_0=a_iE_i$, and
$\nabla_{E_i}E_j=-\delta_{ij}a_iE_0$ give a diagonal curvature operator with
\begin{equation}\label{eq:counterexample-torus-curvature}
 K(E_0,E_i)=-a_i^2,\qquad K(E_i,E_j)=-a_ia_j\quad(i\ne j).
\end{equation}
For a unit vector $E$, $\Ric(E,E)$ is the sum of the sectional curvatures
$K(E,E_j)$ over an orthonormal basis of $E^\perp$. Thus, writing
$\kappa=\alpha^2+(m-1)b^2$,
\begin{align*}
 \Ric(E_0,E_0)&=\sum_{i=1}^m K(E_0,E_i)=-\alpha^2-(m-1)b^2=-\kappa,\\
 \Ric(E_1,E_1)&=K(E_1,E_0)+\sum_{j=2}^m K(E_1,E_j)
              =-\alpha^2-(m-1)\alpha b=-\alpha,\\
 \Ric(E_i,E_i)&=K(E_i,E_0)+K(E_i,E_1)+\sum_{\substack{j\geq2\\j\ne i}}K(E_i,E_j)
              =-b^2-\alpha b-(m-2)b^2=-b\quad(i\geq2),
\end{align*}
where each line uses $\alpha+(m-1)b=\sum_i a_i=1$ to factor
$-\alpha(\alpha+(m-1)b)=-\alpha$ and $-b(\alpha+(m-1)b)=-b$.
The mixed terms vanish by diagonality.
Hence the Ricci eigenspaces $\R E_0$, $\R E_1$,
$\operatorname{span}\{E_2,\ldots,E_m\}$ have eigenvalues
$-\kappa,-\alpha,-b$. Since
\[
 \kappa-b=\alpha^2+(m-1)b^2-b=\alpha(\alpha-b)>0,\qquad
 \alpha-\kappa=\alpha-\alpha^2-(m-1)b^2=(m-1)b(\alpha-b)>0
\]
(using $\alpha+(m-1)b=1$), the smallest eigenvalue is $-\alpha$. Therefore
\begin{equation}\label{eq:counterexample-beta-calculation}
 \rho=-\alpha=-\frac\beta4,\qquad
 \frac{\Delta u}{u}=\frac{u''+u'}u=-\frac14,\qquad
 u^2dV_g=dt\,d\vartheta,
\end{equation}
where $u=e^{-t/2}$ gives $u'/u=-1/2$, $u''/u=1/4$, and
$\Delta u=u''+(\sum_i a_i)u'=u''+u'$. Thus
$L_\beta u=-\beta\Delta u+\rho u=(\beta/4-\alpha)u=0$ since
$\alpha=\beta/4$. The ground-state identity \eqref{eq:ground-state}
gives $Q_\beta\geq0$; the Rayleigh quotient of $u\chi(t/R)$ is
$\beta R^{-2}\int|\chi'|^2/\int\chi^2\to0$ for any nonzero
$\chi\in C_c^\infty(\R)$, proving \eqref{eq:beta-counterexample-zero}.

The three Ricci eigenvalues are distinct. For unit vectors in the three
eigenspaces, respectively, the Jacobi operator has the unequal pairs
\begin{equation*}
 (-\alpha^2,-b^2),\qquad
 (-\alpha^2,-\alpha b),\qquad
 (-b^2,-\alpha b).
\end{equation*}
Indeed, for every unit $V\in\operatorname{span}\{E_2,\ldots,E_m\}$,
diagonality of the curvature operator gives $K(V,E_0)=-b^2$ and
$K(V,E_1)=-\alpha b$. Thus multiplicity in this eigenspace leaves no
additional candidate for a warped radial direction.
The radial-direction criterion therefore excludes every local warped
representation.

Finally, with $\theta=(m-2)/(m-1)$, a fixed torus coordinate gives
\begin{equation}\label{eq:counterexample-incomplete}
 \int_{-\infty}^0u^{-\theta}\,dt=\frac2\theta<\infty.
\end{equation}
Each curve with fixed torus coordinates is a weighted minimizing line:
between two of its points every competitor has weighted length at least
the absolute change of $\int_0^t u(s)^\beta\,ds$, and the radial segment
realizes this bound.
On the complete half-manifold $(-\infty,0]\times\mathbb T^m$, the inward
normal is $\nu=-\partial_t$ and $H=-1$, so
\begin{equation}\label{eq:counterexample-boundary-sign}
 H+\beta\nu(\log u)=-1+\frac\beta2
 =-\frac{m-2}{2(m-1)}<0.
\end{equation}
On this half-manifold $u$ still satisfies $L_\beta u=0$ in the interior
and $q_B<0$ on the boundary. The radial ray toward $t=-\infty$ has
finite $\gstar$-length, and the same curvature calculation excludes
a local warped-product representation with one-dimensional base.
\end{proof}

The hypothesis of Theorem~\ref{thm:closed} allows a choice of positive
supersolution. On the two-ended torus model the ground state determines
that choice up to a constant. Recall that a Schr\"odinger operator $L$
with nonnegative quadratic form is critical if subtracting any nonzero
nonnegative smooth compactly supported potential destroys nonnegativity.
If some such potential can be subtracted while preserving nonnegativity,
$L$ is subcritical.

\begin{corollary}[Positive supersolutions on the torus model]
\label{cor:torus-critical}
On the boundaryless manifold in \eqref{eq:counterexample-torus-metric},
$L_\beta$ is critical. Every global positive weak supersolution
$w\in W^{1,2}_{\mathrm{loc}}(M)$ of $L_\beta w\geq0$ equals $Ce^{-t/2}$
almost everywhere for some $C>0$. With this smooth representative of
$w$, the metric $w^{-2\theta}g$ is incomplete.
\end{corollary}

\begin{proof}
Choose a nonnegative $\chi\in C_c^\infty(\R)$ equal to one on $[-1,1]$.
Write $u=e^{-t/2}$ and $\varphi_R=u\chi(t/R)$. The ground-state identity
and \eqref{eq:counterexample-beta-calculation} give the unnormalized energy
$$
 Q_\beta(\varphi_R)
 =\frac{\beta\,\vol(\mathbb T^m)}{R}
   \int_{\R}|\chi'(s)|^2\,ds .
$$
On every compact set $\varphi_R=u$ for all sufficiently large $R$.
Thus any nonzero $W\in C_c^\infty(M)$ with $W\geq0$ would satisfy
$Q_\beta(\varphi_R)<\int_M W\varphi_R^2\,dV_g$ for large $R$.
This proves criticality by the ground-state alternative applied to
$L_\beta/\beta$ \cite[Theorem~2.3]{CMMR25}.

To see uniqueness directly, let $w$ be a positive weak supersolution.
Testing its inequality with $\varphi^2/(w+\varepsilon)$ and letting
$\varepsilon\downarrow0$ gives, for $\varphi\in C_c^\infty(M)$,
$$
 \beta\int_M\left|\nabla\varphi-\frac{\varphi}{w}\nabla w\right|^2\,dV_g
 \leq Q_\beta(\varphi).
$$
On any precompact open set, $\varphi_R=u$ for all sufficiently
large $R$. Since $Q_\beta(\varphi_R)\to0$, the integrand then gives
$w\nabla u=u\nabla w$ almost everywhere there. The locally Sobolev
function $w/u$ therefore has zero weak gradient. Connectedness gives
$w=Cu$ almost everywhere for some $C>0$. With the representative $Cu$,
the ray in \eqref{eq:counterexample-incomplete} has conformal length
$2C^{-\theta}/\theta$.
\end{proof}

\begin{remark}[What non-removability means]
\label{rem:completeness-nonremovable}
The condition is not necessary for each individual warped product.
Take $g=dt^2+a(t)^2h_{S^m}$ with $a=\sech t$ and
$u=\cosh^{m-1}t$. Since $(a')^2-aa''=\sech^4t\leq1$,
Proposition~\ref{prop:converse} gives $L_\beta u=0$, while each end has
finite $g_*$-length $\int_0^\infty\cosh^{-(m-2)}t\,dt$ for $m>2$.
\end{remark}

\subsection{The critical coefficient}

\begin{proposition}[Critical counterexamples]\label{prop:critical-all}
For every $n\geq4$, there is a complete two-ended metric $g$ on
$\mathbb R\times S^{n-1}$ satisfying
\begin{equation}\label{eq:critical-all-spectrum}
 \lambda_1\left(-\frac4{n-1}\Delta+\rho\right)=0,
\end{equation}
which admits no one-dimensional-base warped-product representation near
some point. For the same metric, $\lambda_1(-\gamma\Delta+\rho)<0$
whenever $0<\gamma<4/(n-1)$.
\end{proposition}

\begin{proof}
Let $m=n-1\geq3$ and let $h_0$ be the unit round metric on $S^m$.
Choose a smooth trace-free $h_0$-selfadjoint endomorphism $A$ with
$\operatorname{tr}(A^2)\leq2$, $A(p)=\operatorname{diag}(1,-1,0,\ldots,0)$,
and $\nabla A(p)=0$. Such an $A$ is obtained by cutting off
$F_1\otimes F_1^\flat-F_2\otimes F_2^\flat$ in a local orthonormal frame
with $\nabla F_i(p)=0$ and extending by zero.
Set $h_z=h_0(e^{2zA}\cdot,\cdot)$. Since $\det(e^{2zA})=1$, $dV_{h_z}=dV_{h_0}$,
and $A$ remains selfadjoint with unchanged eigenvalues.
For sufficiently small $z_*>0$, compactness gives
$\Ric_{h_z}\geq h_z$ and
$|\operatorname{div}_{h_z}A|_{h_z}\leq D$ for $|z|\leq z_*$,
where divergence is the covector $\nabla_i A^i{}_j$.
Choose $C$ with $2(m-2)C/m\geq3+2D^2$, which is possible since $m>2$.
Take $f(t)=\delta t\chi(t)$, with $\chi\in C_c^\infty(\mathbb R)$
equal to one near zero and $\delta>0$ small enough that $|f|\leq z_*$.
Writing $s=f'$, choose
$v_->\frac m2(1+\sqrt{2C/m})\|s\|_\infty$. For $\varepsilon>0$ to be chosen, set
\begin{gather}\label{eq:critical-all-functions}
 v(t)=v_-+C\int_{-\infty}^t s(\tau)^2\,d\tau,\qquad
 u(t)=\exp\left(\int_0^t v(\tau)\,d\tau\right),\qquad a=\varepsilon u^{-2/m},\\
 g=dt^2+a^2h_f,\qquad k=-2v/m,\qquad
 S=k\operatorname{Id}+sA,\qquad H=-2v.\nonumber
\end{gather}
Here $S,H$ are the slice shape operator and mean curvature.
Put $P=(4/m)\Delta u/u=(4/m)(v'-v^2)$, with index $0$ denoting $\partial_t$.
Since $[A,S]=0$, $A$ is parallel along $\partial_t$. The formulas
$\Ric_{00}=-H'-|S|^2$, contracted Codazzi, and
$\Ric_{\rm tan}^{\sharp}=a^{-2}\Ric_{h_f}^{\sharp}-S'-HS$ yield
\begin{align*}
 \Ric_{00}-P
   &=\left(\tfrac{2(m-2)}mC-\operatorname{tr}(A^2)\right)s^2,\\
 |\Ric(\partial_t,\cdot)|_{\rm tan}&\leq a^{-1}D|s|,\\
 T:=(\Ric-Pg)_{\rm tan}^{\sharp}
   &=a^{-2}\Ric_{h_f}^{\sharp}-\tfrac{2C}{m}s^2\operatorname{Id}
       +(2vs-s')A.
\end{align*}
The mixed estimate uses $d_{g_t}H=0$ and
$\operatorname{div}_{g_t}S=s\operatorname{div}_{h_f}A$.
Write $E=T-a^{-2}\Ric_{h_f}^{\sharp}$, and let $I$ be a compact
interval containing the supports of $s,s'$. The operator norm of $E$
is the same for $a^2h_f$ and $h_f$; neither $E$ nor $u$ depends on
$\varepsilon$. Thus $K=\sup_{I\times S^m}u^{-4/m}\|E\|$ is finite.
Choose $\varepsilon^2 K\leq1/2$. Then $a^2\|E\|\leq1/2$ on $I$,
and $E=0$ off $I$, giving $T\geq(2a^2)^{-1}\operatorname{Id}$ everywhere.
Its Schur complement is at least
\begin{equation*}
 \left(\tfrac{2(m-2)}mC-2-2D^2\right)s^2\geq s^2.
\end{equation*}
Thus $\rho\geq(4/m)\Delta u/u$, and \eqref{eq:ground-state} gives
$\lambda_1(-(4/m)\Delta+\rho)\geq0$.
Sufficiently far along the negative end we have $s=s'=v'=0$ and $v=v_-$.
There the mixed Ricci terms vanish, $\Ric_{00}=P$, and the tangential
Ricci eigenvalues are strictly greater than $P$ by the bound on $T$.
Hence $\rho=P=-(4/m)v_-^2$ there. Since $u^2a^m=\varepsilon^m$,
cutoffs $u\chi_R$ supported in this region have Rayleigh quotient
\begin{equation}\label{eq:critical-all-rayleigh}
 \left(\gamma-\frac4m\right)v_-^2
 +\gamma\frac{\int|\chi_R'|^2}{\int\chi_R^2},
 \qquad \frac{\int|\chi_R'|^2}{\int\chi_R^2}\longrightarrow0.
\end{equation}
For example take $\chi_R(t)=\psi((t+3R)/R)$ with
$0\ne\psi\in C_c^\infty((-1,1))$.
Its support lies in $(-4R,-2R)$, where $v=v_-$ for $R$ large.
The cross term vanishes since $\int\chi_R\chi_R'=0$, and the quotient
of cutoff integrals is $R^{-2}\int|\psi'|^2/\int\psi^2$.
This proves both spectral assertions.

At $(0,p)$, $f=0$, $s=\delta$, $s'=0$, and the mixed Ricci components
vanish by $\nabla A(p)=0$. The radial Jacobi operator is
\begin{equation}\label{eq:critical-all-jacobi}
 J_{\partial_t}=\frac{2C}{m}\delta^2\operatorname{Id}
              -(k\operatorname{Id}+\delta A(p))^2.
\end{equation}
At this point each eigenvalue is $(2C/m)\delta^2-(k+\delta\lambda)^2$
with $\lambda\in\{-1,0,1\}$. Since $\delta=s(0)\leq\|s\|_\infty$
and $v(0)\geq v_-$, we have $|k+\delta\lambda|\geq2v_-/m-\delta>
\sqrt{2C/m}\,\delta$. Thus every eigenvalue is negative; the two
associated with $\pm1$ differ by $-4k\delta>0$.
Shrinking $\varepsilon$ further makes the tangential Ricci eigenvalues
larger than the radial one and gives
$K_g(X,Y)\geq\varepsilon^{-2}-\|S\|_{\mathrm{op}}^2>0$
for orthonormal tangent vectors $X,Y$.
Hence every Ricci eigendirection is radial or slice-tangent.
The radial Jacobi operator is nonscalar, while that of any tangent
direction takes a negative value on $\partial_t$ and a positive value
on a perpendicular tangent direction. The criterion above excludes
every possible warped-product radial direction.
\end{proof}

The same spherical model addresses a question about the potential in a
spectral Ricci bound. After Example~3.8 of \cite{CMMR25},
Catino, Mari, Mastrolia, and Roncoroni raise the possibility that on a
manifold with at least two ends the bounds $\Ric_g\geq-cVg$ and
$-\Delta-V\geq0$ at $c=4/(n-1)$ already force criticality and the
identity $\rho+cV=0$ everywhere. Their Theorem~3.3 gives this conclusion
in dimension three. In dimension $n\geq4$ it assumes a positive
supersolution $w$ of $(-\Delta-V)w\geq0$ outside a compact set with
$w(x)=o(r(x)^{(n-1)^2/(4(n-3))})$, where $r$ is distance from a fixed
point. The strict Ricci inequality on the deformed slices allows a
compactly supported perturbation of the potential.

\begin{corollary}[A subcritical operator with zero spectral bottom]
\label{cor:sphere-subcritical}
For every $n\geq4$, the metric of Proposition~\ref{prop:critical-all}
admits a smooth potential $V$ for which
$$
 \Ric_g\geq-cVg,\qquad
 \lambda_1(-\Delta-V)=0,\qquad c=\frac{4}{n-1}.
$$
The operator $-\Delta-V$ is subcritical, and $\rho+cV>0$ on a
nonempty open set.
\end{corollary}

\begin{proof}
Retain $u$, $s=f'$ and $P=c\Delta u/u$ from the preceding construction.
The estimate for the tangential block of $\Ric_g-Pg$ gives a positive
definite block everywhere, and its Schur complement is bounded below
by $s^2$. Hence $G=\rho-P$ is continuous and nonnegative, and it is
strictly positive on $U=\{s\ne0\}$, which contains the slice $t=0$.
Choose a nonzero $W\in C_c^\infty(U)$ with $W\geq0$ and scale it so
that $cW\leq G/2$. Set $V_0=-\Delta u/u$, $V=V_0-W$, and $L=-\Delta-V$.
Then $\rho+cV=G-cW\geq G/2$. Since $(-\Delta-V_0)u=0$, we have
$$
 Q_L(\varphi)-\int_M W\varphi^2\,dV_g
 =\int_M u^2\left|\nabla(\varphi/u)\right|^2\,dV_g\geq0
 \qquad\bigl(\varphi\in C_c^\infty(M)\bigr),
$$
where $Q_L$ is the quadratic form of $L$. Thus $W$ is a nonzero Hardy
weight for $L$, proving subcriticality.

It remains to locate the spectral bottom. For
$0\ne\psi\in C_c^\infty((-1,1))$, the cutoff $\chi_R(t)=\psi((t+3R)/R)$
has support disjoint from $W$ for large $R$. Since
$u^2dV_g=\varepsilon^m\,dt\,dV_{h_0}$, the preceding identity yields
$$
 \frac{Q_L(u\chi_R)}{\int_M u^2\chi_R^2\,dV_g}
 =\frac{1}{R^2}
   \frac{\int_{\R}|\psi'|^2\,dt}{\int_{\R}\psi^2\,dt}
 \longrightarrow 0.
$$
Together with $L\geq0$ this proves the claim.
\end{proof}

\section{Scalar-flat three-manifolds and spatial Schwarzschild}\label{sec:applications}

This section applies the splitting theorems to scalar-flat
three-manifolds. We first show that scalar flatness gives a first
integral for the warping function in the splitting coordinates, and that
the radial Ricci equality makes the associated constant nonnegative.
Normalizing the spherical fiber to unit curvature then identifies the
Schwarzschild mass. Combining these ingredients proves the two
classifications of this section: flat torus products and positive-mass
spatial Schwarzschild doubles in the two-ended case, and flat
half-cylinders, Euclidean ball exteriors, and regions of Schwarzschild
doubles bounded by coordinate spheres in the compact boundary case.

\begin{proposition}[Two-ended scalar flat manifolds]\label{prop:schwarz-closed}
Let $(M^3,g)$ be smooth, complete, connected, orientable, and without
boundary. Suppose $R_g=0$ and $M$ has at least two ends. Then
$Q_2\geq0$ if and only if $M$ is isometric either to the product of
$\R$ with a flat 2-torus or to the doubled spatial Schwarzschild
manifold of mass $\mathfrak m>0$. The latter is the smooth double
across $r=2\mathfrak m$ of the metric
$$
 g_{\mathfrak m}=\frac{dr^2}{1-2\mathfrak m/r}
                  +r^2g_{S^2(1)},\qquad r>2\mathfrak m .
$$
\end{proposition}

\begin{proof}
By Corollary~\ref{cor:three}, $g=dt^2+a(t)^2h$ on $\R\times\Sigma$,
where $\Sigma$ is a closed connected orientable surface. The radial
equality \eqref{eq:radial-equality} of Theorem~\ref{thm:closed} gives
$\rho=\Ric(\partial_t,\partial_t)$. For any unit vector $X$ tangent to a
slice, the curvature formulas give
$$
 \Ric(\partial_t,\partial_t)=-2\frac{a''}{a},\qquad
 \Ric(X,X)=\frac{K_h-(a')^2-aa''}{a^2}.
$$
Taking the trace and using $R_g=0$ gives
$K_h(x)=(a'(t))^2+2a(t)a''(t)$. Both sides equal a constant $k$.
Substitution in the Ricci difference yields
$\Ric(X,X)-\Ric(\partial_t,\partial_t)=3a''/a\geq0$. Consequently
$$
 \mu=\frac{a}{2}\big(k-(a')^2\big)=a^2a''\geq0,\qquad
 \mu'=\frac{a'}{2}\big(k-(a')^2-2aa''\big)=0 .
$$
If $\mu=0$, the positive affine function $a$ is constant on $\R$, and
$k=0$. Thus the fiber is a flat two-torus. If $\mu>0$, then $k>0$
and $(\Sigma,kh)$ is a round unit sphere. Its area radius $r=a/\sqrt{k}$
and the constant $\mathfrak m=\mu/k^{3/2}$ satisfy
$$
 (r')^2=1-\frac{2\mathfrak m}{r},\qquad
 r''=\frac{\mathfrak m}{r^2},\qquad r\geq2\mathfrak m .
$$
The function $r$ is strictly convex and $|r'|\leq1$. It tends to
infinity at both ends, since boundedness on a half-line would give a
positive lower bound for $r''$ there. At its unique minimum it equals
$2\mathfrak m$. Uniqueness for the second-order equation makes $r$ even
about this point. On either side $r$ is a coordinate and gives the
metric $g_{\mathfrak m}$ in the statement. In the coordinate
$r=s(1+\mathfrak m/(2s))^2$, $s>\mathfrak m/2$, the metric becomes
$(1+\mathfrak m/(2s))^4(ds^2+s^2g_{S^2(1)})$. Thus both ends are
asymptotically flat and have ADM mass $\mathfrak m$.

Conversely, the product is flat. On the double, $\rho=-2r''/r$ and
$u=1/r$ satisfies $\Delta u=-r''/r^2$ everywhere, including the
minimal sphere. Hence $L_2u=0$, and the ground-state identity
\eqref{eq:ground-state} gives $Q_2\geq0$. This is the warped model
calculation of \cite[Proposition~6.1]{Li-Wang-ASENS-weighted} and
\cite[Example~3.7]{CMMR25}.
\end{proof}

With compact boundary the same constant determines the mass, while the
sign of the inward mean curvature specifies which side of a coordinate
sphere occurs. The spectral hypothesis then includes the boundary term
arising from integration by parts.

\begin{proposition}[Scalar flat manifolds with compact boundary]
\label{prop:schwarz-boundary}
Let $(M^3,g)$ be a smooth, complete, connected, orientable, noncompact
manifold with $R_g=0$ and with nonempty compact connected boundary
$B=\partial M$. Write $H$ for the sum of the principal curvatures
with respect to the inward unit normal $\nu$. Then
$$
 Q_B(\varphi):=\int_M(2|\nabla\varphi|^2+\rho\varphi^2)\,dV_g
                    -\int_B H\varphi^2\,dA_g\geq0
$$
for all smooth, compactly supported $\varphi$ up to $B$, with arbitrary
boundary trace, if and only if $M$ is one of the following:
\begin{enumerate}
\item the product of $[0,\infty)$ with a flat two-torus,
\item the closed exterior of a round ball in Euclidean three-space,
\item the closure of either component obtained by cutting a positive-mass
      doubled spatial Schwarzschild manifold along a coordinate sphere,
      including the minimal sphere.
\end{enumerate}
If $B$ is minimal, the possibilities reduce to the flat half-cylinder
and the Schwarzschild exterior bounded by its horizon. In particular,
a minimal spherical boundary determines the latter model with
$\mathfrak m=\sqrt{|B|/(16\pi)}$.
\end{proposition}

\begin{proof}
We first construct a supersolution to which Theorem~\ref{thm:boundary}
applies. Exhaust $M$ by smooth, precompact connected domains
containing a collar of $B$. On each domain impose Dirichlet data on the
artificial boundary and the natural Robin condition on $B$. By zero
extension across the artificial boundary and density in the form
domain, the first eigenvalues satisfy $\lambda_j\geq0$. Domain
monotonicity gives $\lambda_j\downarrow\lambda\geq0$. Normalize the
positive first eigenfunctions to equal one at a fixed interior point.
Interior and boundary Harnack estimates give local bounds uniform in $j$.
For the boundary estimate one can multiply the eigenfunction by a
smooth positive factor whose inward logarithmic derivative is $H/2$
on $B$. The resulting homogeneous Neumann condition permits even
reflection in boundary normal coordinates and application of the
interior Harnack inequality. Local elliptic estimates for the Robin
problem now give a subsequence converging in $C^2$ on compact subsets
up to $B$. Its limit $u$ is positive up to $B$, is locally
$C^{2,\alpha}$ for $0<\alpha<1$, and satisfies
$$
 L_2u=\lambda u\quad\text{in }\interior,\qquad
 2\partial_{\rm out}u-Hu=0\quad\text{on }B .
$$
Thus $H+2\nu(\log u)=0$. Conversely, for any positive $u$ smooth
enough for integration by parts one has
$$
 \begin{aligned}
 Q_B(\varphi)={}&2\int_M u^2|\nabla(\varphi/u)|^2\,dV_g
       +\int_M\frac{L_2u}{u}\varphi^2\,dV_g\\
       &-\int_B\bigl(H+2\nu(\log u)\bigr)\varphi^2\,dA_g .
 \end{aligned}
$$
In particular a positive interior supersolution with
$H+2\nu(\log u)\leq0$ gives $Q_B\geq0$.

Theorem~\ref{thm:boundary} yields $g=dt^2+a(t)^2h$ for $t\geq0$,
with $h=g|_B$, $a(0)=1$, and radial Ricci curvature equal to $\rho$.
The preceding calculation of $k$ and $\mu$ applies on this half-line.
If $\mu=0$, then $a(t)=1+bt$ with $b\geq0$ and $k=b^2$. The case $b=0$ gives the flat
half-cylinder. For $b>0$, $(B,b^2h)$ is the unit sphere and $r=a/b$
gives the Euclidean exterior. If $\mu>0$, the same normalization gives
$r''=\mathfrak m/r^2$ and $(r')^2=1-2\mathfrak m/r$. The doubled
solution from Proposition~\ref{prop:schwarz-closed} realizes every
allowed pair $(r(0),r'(0))$. Uniqueness for this initial-value problem
identifies the half-line with the corresponding side of a coordinate
sphere. Each listed model has $u=1/a$, for which $L_2u=0$ and
$H+2\nu(\log u)=0$, proving the converse.

For a spherical boundary set $r_B=\sqrt{|B|/(4\pi)}$. Then
$H=2r'(0)/r_B$ is constant, and its area and mean curvature give the
mass, including the Euclidean value zero, by
$$
 \mathfrak m=\frac{r_B}{2}\left(1-\frac{H^2r_B^2}{4}\right).
$$
For $H>0$ the radius increases from the boundary. For $H<0$ it first
decreases to the minimal sphere and then increases toward infinity.
The case $H=0$ in the spherical branch starts at the horizon and gives
the asserted mass--area identity.
\end{proof}

The spectral conditions above determine spherical symmetry and the
asymptotically flat geometry of the Schwarzschild branch. Other
characterizations start from different geometric data. Reiris
Ithurralde~\cite{ReirisI,ReirisII} classifies complete static vacuum
black holes with compact horizon and derives the possible end
geometries. On $\R^3\setminus\{0\}$, Miao assumes nonnegative scalar
curvature, asymptotic flatness at infinity, and a harmonic function
tending to zero at the puncture and to one at infinity. Under these
hypotheses equality in the mass--capacity estimate obtained as the
level tends to zero characterizes the entire doubled Schwarzschild
manifold \cite[Theorem~1.5(ii)]{MiaoHarmonic}.
The mass--area identity above follows
from the identified metric. The Riemannian Penrose inequality of
Huisken and Ilmanen~\cite{HuiskenIlmanen} gives a mass lower bound with
the same constant for asymptotically flat three-manifolds of nonnegative
scalar curvature with connected outermost minimal boundary.

\bibliographystyle{alpha}
\bibliography{splitingreference}
\end{document}